\documentclass{article}
\usepackage{amssymb}
\usepackage{amsfonts}
\usepackage{amsmath}

\input{tcilatex}
\begin{document}

\title{On Bernstein processes of maximal entropy }
\author{Pierre-A. Vuillermot$^{\ast ,\ast \ast }$ \\
UMR-CNRS 7502, Inst. \'{E}lie Cartan de Lorraine, Nancy, France$^{\ast }$\\
Grupo de F\'{\i}sica Matem\'{a}tica, GFMUL, Faculdade de Ci\^{e}ncias, \\
Universidade de Lisboa, 1749-016 Lisboa, Portugal$^{\ast \ast }$}
\date{}
\maketitle

\begin{abstract}
In this article we define and investigate statistical operators and an
entropy functional for Bernstein stochastic processes associated with
hierarchies of forward-backward systems of decoupled deterministic linear
parabolic partial differential equations. The systems under consideration
are defined on open bounded domains $D\subset \mathbb{R}^{d}$ of Euclidean
space where $d\in \mathbb{N}^{+}$ is arbitrary, and are subject to Neumann
boundary conditions. We assume that the elliptic part of the parabolic
operator in the equations is a self-adjoint Schr\"{o}dinger operator,
bounded from below and with compact resolvent in $L^{2}(D)$. The statistical
operators we consider are then trace-class operators defined from sequences
of probabilities associated with the point spectrum of the elliptic part in
question, which allow the distinction between pure and mixed processes. We
prove in particular that the Bernstein processes of maximal entropy are
those for which the associated sequences of probabilities are of Gibbs type.
We illustrate our results by considering processes associated with a
specific hierarchy of forward-backward heat equations defined in a
two-dimensional disk.

\ \ \ \ \ \ \ \ \ \ \ \ \ \ \ \ \ \ \ \ \ \ \ \ \ \ \ \ \ \ \ \ \ \ \ \ \ \
\end{abstract}

\section{Introduction and outline}

The theory of Bernstein (or reciprocal) processes was launched many years
ago in \cite{bernstein} following the seminal contribution put forward in
the last section of \cite{schroedinger}. At the very end of \cite%
{schroedinger}, Schr\"{o}dinger indeed gave a positive answer to the
question whether it is possible to generate a reversible diffusion process
from a pair of adjoint, deterministic, linear parabolic partial differential
equations whose solutions typically display irreversible behavior. The
considerations of \cite{schroedinger} were based on entropy arguments, and
have had many important ramifications and generalizations over the years up
to this day, including connections with Optimal Transport Theory (see, e.g., 
\cite{beurling}-\cite{villani}, and the many references therein). On the
other hand, a systematic and abstract investigation of continuous time
versions of the processes was carried out in \cite{jamison}, according to
which it became clear that Bernstein processes may exist without any
reference to partial differential equations and may admit as state space any
topological space countable at infinity. In spite of that, a great deal of
attention has recently been paid to the way that such processes may be
generated in Euclidean space of arbitrary dimension from certain particular
systems of parabolic partial differential equations, thereby allowing one to
recast things within the original framework of \cite{schroedinger} with the
goal of investigating those processes that are not Markovian (see, e.g., 
\cite{vuillermot}-\cite{vuillerzamb}).

It is our purpose here to continue and deepen our analysis of such
processes, and accordingly we shall organize the remaining part of this
article in the following way: in Section 2 we introduce a hierarchy of
forward-backward systems of decoupled, deterministic, linear parabolic
partial differential equations defined on open bounded domains of Euclidean
space. Those systems are characterized by the fact that the elliptic part of
the parabolic operator is, up to a sign, a self-adjoint Schr\"{o}dinger
operator bounded from below and with compact resolvent in standard $L^{2}$%
-space. The hierarchy comes about by associating with each level of the pure
point spectrum of the elliptic part a suitable pair of initial-final data.
We then proceed by defining what a Bernstein process is, and show how we can
construct from the hierarchy we just alluded to a sequence of such processes
that are Markovian. This requires the existence of probability measures of a
very specific form which we obtain from the initial-final data and the heat
kernel of the given system. In Section 2 we also associate with the spectrum
of the elliptic part a sequence of probabilities which eventually allows us
to construct non-Markovian processes by means of a suitable averaging
procedure, as well as the related statistical operators and the entropy
functional which we investigate in detail. Those operators are important in
that they allow the classification of the processes as pure or mixed, and we
prove in particular that the Bernstein processes of maximal entropy are
those for which the probabilities in question are of Gibbs type. In Section
3 we illustrate some of our results by considering Bernstein processes
generated by a specific hierarchy of forward-backward heat equations and
wandering in a two-dimensional disk, ending up with fairly explicit formulae
for the corresponding probabilities and expectation values. Finally, we
devote Appendix A to the analysis of statistical operators which are more
general than that investigated in Section 2, and Appendix B to stating a
general result regarding the very existence of Bernstein processes that goes
back to \cite{jamison} and \cite{vuillerzambis}, which we slightly
reformulated for the needs of this article. We conclude Appendix B with a
brief remark regarding the connection between Bernstein processes, Schr\"{o}%
dinger's problem and Optimal Transport Theory.

\section{Statistical operators and an entropy functional for Bernstein
processes}

Let $D\subset \mathbb{R}^{d}$ with $d\in \mathbb{N}^{+}$be an open bounded
domain with a sufficiently smooth boundary $\partial D$ and let $L^{2}(D)$
be the standard Hilbert space of all Lebesgue-measurable, square-integrable
complex-valued functions on $D$ with respect to Lebesgue measure, whose
inner product and induced norm we shall denote by $(.,.)_{2}$ and $%
\left\Vert .\right\Vert _{2}$, respectively. Let us consider the
differential operator%
\begin{equation}
\mathcal{H}=-\frac{1}{2}\Delta _{\mathsf{x}}+V  \label{hamiltonian}
\end{equation}%
where $\Delta _{\mathsf{x}}$ stands for Neumann's Laplacian on $D$ and where
the following hypothesis holds for the additional term:

\bigskip

(H$_{1}$) The function $V:D\mapsto \mathbb{R}$ satisfies $V\in L^{p}(D)$
where

\bigskip

\begin{equation*}
p\in \left\{ 
\begin{array}{c}
\left[ 1,+\infty \right] \text{ \ \ if }d=1, \\ 
\\ 
\left( 1,+\infty \right] \text{ \ \ if }d=2, \\ 
\\ 
\left[ \frac{d}{2},+\infty \right] \text{ if }d\geq 3%
\end{array}%
\right.
\end{equation*}%
and is bounded from below.

\bigskip

Under these conditions it is well known that (\ref{hamiltonian}) admits a
self-adjoint realization with compact resolvent in $L^{2}(D)$ and thereby a
pure point spectrum $\left( \lambda _{\mathsf{m}}\right) _{\mathsf{m}\in 
\mathbb{N}^{d}}$ such that $\lambda _{\mathsf{m}}\rightarrow +\infty $ as $%
\left\vert \mathsf{m}\right\vert :=\sum_{j=1}^{d}m_{j}\rightarrow +\infty $,
whose corresponding eigenfunctions $\left( \mathsf{f}_{\mathsf{m}}\right) _{%
\mathsf{m}\in \mathbb{N}^{d}}$ constitute an orthonormal basis of $L^{2}(D)$
and are assumed to be real (see, e.g., Chapter VI in \cite{edmundsevans},
particularly Theorem 1.9). For each $\mathsf{m}\in \mathbb{N}^{d}$ and $T\in
\left( 0,+\infty \right) $ arbitrary, we then introduce the system of
adjoint, deterministic, linear parabolic partial differential equations
given by

\begin{align}
\partial _{t}u(\mathsf{x},t)& =\frac{1}{2}\Delta _{\mathsf{x}}u(\mathsf{x}%
,t)-V(\mathsf{x})u(\mathsf{x},t),\text{ \ }(\mathsf{x},t)\in D\times \left(
0,T\right] ,  \notag \\
u(\mathsf{x},0)& =\varphi _{0,\mathsf{m}}(\mathsf{x}),\text{ \ \ }\mathsf{x}%
\in D,  \label{forwardproblem} \\
\frac{\partial u(\mathsf{x},t)}{\partial n\mathsf{(x)}}& =0,\text{ \ \ \ }(%
\mathsf{x},t)\in \partial D\times \left( 0,T\right]  \notag
\end{align}%
and%
\begin{align}
-\partial _{t}v(\mathsf{x},t)& =\frac{1}{2}\Delta _{\mathsf{x}}v(\mathsf{x}%
,t)-V(\mathsf{x})v(\mathsf{x},t),\text{ }(\mathsf{x},t)\in D\times \left[
0,T\right) ,  \notag \\
v(\mathsf{x,}T)& =\psi _{T,\mathsf{m}}(\mathsf{x}),\text{ \ \ }\mathsf{x}\in
D,  \label{backwardproblem} \\
\frac{\partial v(\mathsf{x},t)}{\partial n\mathsf{(x)}}& =0,\text{ \ \ \ }(%
\mathsf{x},t)\in \partial D\times \left[ 0,T\right) ,  \notag
\end{align}%
respectively, where $n\mathsf{(x)}$ stands for the unit outer normal to $%
\partial D$ at the point $\mathsf{x}$ and where $\varphi _{0,\mathsf{m}}$, $%
\psi _{T,\mathsf{m}}$ are real-valued functions to be specified below. In
this way we are thus considering a hierarchy of problems of the form (\ref%
{forwardproblem})-(\ref{backwardproblem}), that is, an infinite sequence of
pairs of such equations where each pair is associated with a level of the
spectrum of (\ref{hamiltonian}) through the initial-final data. Furthermore,
an essential ingredient in the forthcoming considerations will be the heat
kernel (or fundamental solution) associated with (\ref{forwardproblem})-(\ref%
{backwardproblem}), which satisfies%
\begin{equation}
\left\{ 
\begin{array}{c}
g(\mathsf{x},t,\mathsf{y})=g(\mathsf{y},t,\mathsf{x}), \\ 
\\ 
0<g(\mathsf{x},t,\mathsf{y})\leq c_{1}t^{-\frac{d}{2}}\exp \left[ -c_{2}%
\frac{\left\vert \mathsf{x}-\mathsf{y}\right\vert ^{2}}{t}\right]%
\end{array}%
\right.  \label{heatkernel}
\end{equation}%
for all $\mathsf{x},\mathsf{y\in }\mathbb{\ }\overline{D}\mathbb{\ }$and
every $t\in \left( 0,T\right] $ for some $c_{1,2}>0$. It is indeed the
knowledge of $\varphi _{0,\mathsf{m}}$, $\psi _{T,\mathsf{m}}$ and (\ref%
{heatkernel}) that will allow us to construct sequences of Bernstein
processes $Z_{\tau \in \left[ 0,T\right] }^{\mathsf{m}}$ wandering in $%
\overline{D}$. We begin with the following:

\bigskip

\textbf{Definition 1. }\textit{We say the }$\overline{D}$\textit{-valued
process }$Z_{\tau \in \left[ 0,T\right] }$\textit{\ defined on the complete
probability space }$\left( \Omega ,\mathcal{F},\mathbb{P}\right) $\textit{\
is a Bernstein process if }%
\begin{equation}
\mathbb{E}\left( b(Z_{r})\left\vert \mathcal{F}_{s}^{+}\vee \mathcal{F}%
_{t}^{-}\right. \right) =\mathbb{E}\left( b(Z_{r})\left\vert
Z_{s},Z_{t}\right. \right)  \label{conditionalexpectations}
\end{equation}%
$\mathbb{P}$\textit{-almost everywhere for every bounded Borel measurable
function }$b:\overline{D}\mapsto \mathbb{C}$\textit{, and for all }$r,s,t$%
\textit{\ satisfying }$r\in \left( s,t\right) \subset \left[ 0,T\right] $%
\textit{, where }$\mathbb{E}\left( .\left\vert .\right. \right) $\textit{\
denotes the conditional expectation on }$\left( \Omega ,\mathcal{F},\mathbb{P%
}\right) $. \textit{The }$\sigma $\textit{-algebras in (\ref%
{conditionalexpectations}) are}%
\begin{equation*}
\mathcal{F}_{s}^{+}:=\sigma \left\{ Z_{\tau }^{-1}\left( F\right) :\tau \leq
s,\text{ }F\in \mathcal{B}(\overline{D})\right\}
\end{equation*}%
\textit{and }%
\begin{equation*}
\mathcal{F}_{t}^{-}:=\sigma \left\{ Z_{\tau }^{-1}\left( F\right) :\tau \geq
t,\text{ }F\in \mathcal{B}(\overline{D})\right\} ,
\end{equation*}%
\textit{respectively, where }$\mathcal{B}(\overline{D})$ \textit{stands for} 
\textit{the Borel }$\sigma $\textit{-algebra over }$\overline{D}$\textit{.}

\bigskip

The preceding definition is just one out of many equivalent ways of defining
a Bernstein process (see, e.g. \cite{jamison}). It shows that as soon as $%
Z_{s}$ and $Z_{t}$ are known, the behavior of such a process for $\tau \in %
\left[ s,t\right] $ is independent of the statistical information available
prior to time $s$ and after time $t$ as encoded in $\mathcal{F}_{s}^{+}$ and 
$\mathcal{F}_{t}^{-}$, respectively. In fact, a simple probabilistic
argument implies that Relation (\ref{conditionalexpectations}) is equivalent
to the statement that the $\sigma $-algebra 
\begin{equation*}
\mathcal{F}_{\left[ s,t\right] }:=\sigma \left\{ Z_{\tau }^{-1}\left(
F\right) :\tau \in \left[ s,t\right] ,\text{ }F\in \mathcal{B}(\overline{D}%
)\right\}
\end{equation*}%
is conditionally independent of $\mathcal{F}_{s}^{+}\vee \mathcal{F}_{t}^{-}$
when%
\begin{equation*}
\mathcal{F}_{\left\{ s,t\right\} }:=\sigma \left\{ Z_{s}^{-1}\left( F\right)
,\text{ }Z_{t}^{-1}\left( F\right) \text{ }:F\in \mathcal{B}(\overline{D}%
)\right\}
\end{equation*}%
is given (see, e.g., Section 25 of Chapter VII in \cite{loeve} for the
notion of conditionally independent $\sigma $-algebra). Aside from this
property which generalizes Markov's, it is also clear that the above
definition maintains a perfect symmetry between past and future in that the $%
\sigma $-algebras $\mathcal{F}_{s}^{+}$ and $\mathcal{F}_{t}^{-}$ play an
identical r\^{o}le. Let us now assume that $\varphi _{\mathsf{m,0}}>0$ and $%
\psi _{\mathsf{m,}T}>0$ are sufficiently smooth on $\overline{D}$ and let us
consider the probability measures%
\begin{equation}
\mu _{\mathsf{m}}(G)=\int_{G}\mathsf{dxdy}\varphi _{\mathsf{m,0}}(\mathsf{x)}%
g(\mathsf{x},T,\mathsf{y})\psi _{\mathsf{m,}T}(\mathsf{y})
\label{probabilitymeasures}
\end{equation}%
for every $G\in \mathcal{B}(\overline{D})\times \mathcal{B}(\overline{D})$,
which satisfy%
\begin{equation}
\int_{D\times D}\mathsf{dxdy}\varphi _{\mathsf{m,0}}(\mathsf{x)}g(\mathsf{x}%
,T,\mathsf{y})\psi _{\mathsf{m,}T}(\mathsf{y})=1  \label{normalization}
\end{equation}%
where $g$ is the heat kernel (\ref{heatkernel}) pinned down at the terminal
time $T$. Then, writing 
\begin{equation}
u_{\mathsf{m}}(\mathsf{x},t)=\int_{D}\mathsf{dy}g(\mathsf{x},t,\mathsf{y}%
)\varphi _{\mathsf{m,0}}(\mathsf{y)>0\label{forwardsolution}}
\end{equation}%
for the solution to (\ref{forwardproblem}) and%
\begin{equation}
v_{\mathsf{m}}(\mathsf{x},t)=\int_{D}\mathsf{dy}g(\mathsf{x},T-t,\mathsf{y}%
)\psi _{\mathsf{m,}T}(\mathsf{y)>0\label{backwardsolution}}
\end{equation}%
for the solution to (\ref{backwardproblem}), we have the following result
which follows from the substitution of (\ref{probabilitymeasures}) into the
formulae of Theorem B.1 of Appendix B, and from Theorem 2 in \cite%
{vuillerzambis} as far as the Markov property is concerned:

\bigskip

\textbf{Theorem 1. }\textit{Assume that Hypothesis }(H$_{1}$)\textit{\
holds. Then for every }$\mathsf{m}\in \mathbb{N}^{d}$ \textit{there exists a
probability space }$\left( \Omega ,\mathcal{F},\mathbb{P}_{\mu _{\mathsf{m}%
}}\right) $\textit{\ supporting a }$\overline{D}$-\textit{valued} \textit{%
Bernstein process} $Z_{\tau \in \left[ 0,T\right] }^{\mathsf{m}}$ \textit{%
such that the following statements are valid:}

\textit{(a) The process }$Z_{\tau \in \left[ 0,T\right] }^{\mathsf{m}}$ 
\textit{is a forward Markov process whose finite-dimensional distributions
are}%
\begin{eqnarray}
&&\mathbb{P}_{\mu _{\mathsf{m}}}\left( Z_{t_{1}}^{\mathsf{m}}\in
F_{1},...,Z_{t_{n}}^{\mathsf{m}}\in F_{n}\right)  \notag \\
&=&\int_{D}\mathsf{dx}\rho _{\mathsf{m},0}(\mathsf{x})\int_{F_{1}}\mathsf{dx}%
_{1}...\int_{F_{n}}\mathsf{dx}_{n}\dprod\limits_{k=1}^{n}w_{\mathsf{m}%
}^{\ast }\left( \mathsf{x}_{k-1},t_{k-1};\mathsf{x}_{k},t_{k}\right)
\label{distributionter}
\end{eqnarray}%
\textit{for every }$n\in \mathbb{N}^{+}$,\textit{\ all }$F_{1},...,F_{n}\in $%
\textit{\ }$\mathcal{B}(\overline{D})$ \textit{and all }$%
0=t_{0}<t_{1}<...<t_{n}<T$\textit{, with }$\mathsf{x}_{\mathsf{0}}=\mathsf{x}
$. \textit{In the preceding expression the density of the forward Markov
transition function is}%
\begin{equation}
w_{\mathsf{m}}^{\ast }(\mathsf{x},s;\mathsf{y},t)=g(\mathsf{x},t-s,\mathsf{y}%
)\frac{v_{\mathsf{m}}(\mathsf{y},t)}{v_{\mathsf{m}}(\mathsf{x},s)}
\label{markovdensity}
\end{equation}%
\textit{for all} $\mathsf{x,y}\in $ $\overline{D}$ \textit{and all} $s,t\in %
\left[ 0,T\right] $ \textit{with }$t>s$, \textit{while the initial
distribution of the process reads}%
\begin{equation}
\rho _{\mathsf{m},0}(\mathsf{x})=\varphi _{\mathsf{m,0}}(\mathsf{x})v_{%
\mathsf{m}}(\mathsf{x},0).  \label{initialmarginal}
\end{equation}

\textit{(b) The process }$Z_{\tau \in \left[ 0,T\right] }^{\mathsf{m}}$ 
\textit{may also be viewed as a backward Markov process since the
finite-dimensional distributions (\ref{distributionter}) may also be written
as}%
\begin{eqnarray}
&&\mathbb{P}_{\mu _{\mathsf{m}}}\left( Z_{t_{1}}^{\mathsf{m}}\in
F_{1},...,Z_{t_{n}}^{\mathsf{m}}\in F_{n}\right)  \notag \\
&=&\int_{D}\mathsf{dx}\rho _{\mathsf{m},T}(\mathsf{x})\int_{F_{1}}\mathsf{dx}%
_{1}...\int_{F_{n}}\mathsf{dx}_{n}\dprod\limits_{k=1}^{n}w_{\mathsf{m}%
}\left( \mathsf{x}_{k+1},t_{k+1};\mathsf{x}_{k},t_{k}\right)
\label{distributionquarto}
\end{eqnarray}%
\textit{for every }$n\in \mathbb{N}^{+}$,\textit{\ all }$F_{1},...,F_{n}\in $%
\textit{\ }$\mathcal{B}(\overline{D})$ \textit{and all }$%
0<t_{1}<...<t_{n}<t_{n+1}=T$\textit{, with }$\mathsf{x}_{n\mathsf{+1}}=%
\mathsf{x}$. \textit{In the preceding expression the density of the backward
Markov transition function is}%
\begin{equation}
w_{\mathsf{m}}(\mathsf{x},t;\mathsf{y},s)=g(\mathsf{x},t-s,\mathsf{y})\frac{%
u_{\mathsf{m}}(\mathsf{y},s)}{u_{\mathsf{m}}(\mathsf{x},t)}
\label{markovdensitybis}
\end{equation}%
\textit{for all} $\mathsf{x,y}\in $ $\overline{D}$ \textit{and all} $s,t\in %
\left[ 0,T\right] $ \textit{with }$t>s$, \textit{while the final
distribution of the process reads}%
\begin{equation*}
\rho _{\mathsf{m},T}(\mathsf{x})=\psi _{\mathsf{m,}T}(\mathsf{x})u_{\mathsf{m%
}}(\mathsf{x},T).
\end{equation*}

\textit{(c) We have} 
\begin{equation}
\mathbb{P}_{\mu _{\mathsf{m}}}\left( Z_{t}^{\mathsf{m}}\in F\right) =\int_{F}%
\mathsf{dx}u_{\mathsf{m}}(\mathsf{x},t)v_{\mathsf{m}}(\mathsf{x},t)
\label{absoluteprobability}
\end{equation}%
\textit{for each }$t\in \left[ 0,T\right] $ \textit{and every} $F\in 
\mathcal{B}(\overline{D})$.

\textit{(d) Finally,}%
\begin{equation}
\mathbb{E}_{\mu _{\mathsf{m}}}\left( b(Z_{t}^{\mathsf{m}})\right) =\int_{D}%
\mathsf{dx}b(\mathsf{x})u_{\mathsf{m}}(\mathsf{x},t)v_{\mathsf{m}}(\mathsf{x}%
,t)  \label{expectations}
\end{equation}%
\textit{for each bounded Borel measurable function }$b:\overline{D}$ $%
\mathbb{\mapsto C}$ \textit{and every }$t\in \left[ 0,T\right] $.

\bigskip

\textsc{Remarks. }$\left( \text{1}\right) $\textsc{\ }The fact that $Z_{\tau
\in \left[ 0,T\right] }^{\mathsf{m}}$ is a Markov process for each $\mathsf{m%
}$ may be read off Relations (\ref{distributionter}) and (\ref%
{distributionquarto}), inasmuch as (\ref{markovdensity}) and (\ref%
{markovdensitybis}) are the densities of transition functions that satisfy
the Chapman-Kolmogorov equation (see, e.g., Lemmas 1 and 2 in \cite%
{vuillerzambis}, and for more general comments Section 2.4 in Chapter 2 of 
\cite{ito}). Alternatively, the Markov property of $Z_{\tau \in \left[ 0,T%
\right] }^{\mathsf{m}}$ is an immediate consequence of the form (\ref%
{probabilitymeasures}) of the underlying probability measures through
Theorem 3.1 in \cite{jamison}. Furthermore, the fact that $Z_{\tau \in \left[
0,T\right] }^{\mathsf{m}}$ is both a forward and a backward Markov process
is related to the perfect symmetry between past and future which we alluded
to above, also encoded in (\ref{absoluteprobability}) where (\ref%
{forwardsolution}) and (\ref{backwardsolution}) play an equivalent r\^{o}le.
We refer the reader to \cite{vuillerzambis} for further considerations on
this issue, where a general notion of reversibility was put forward in order
to deal with processes generated by systems of non-autonomous
forward-backward parabolic equations. Finally, we note that the processes $%
Z_{\tau \in \left[ 0,T\right] }^{\mathsf{m}}$ are in general non-stationary
(see, e.g., our construction in Section 3).

$\left( \text{2}\right) $ Theorem 3.1 in \cite{jamison} actually says much
more than what we just referred to in the preceding remark. Indeed, when
applied to the present situation, it asserts that one may generate a
Markovian Bernstein process from a probability measure $\mu $ on $\mathcal{B}%
(\overline{D})\times \mathcal{B}(\overline{D})$ if, and only if, there exist
positive measures $\nu _{0}$ and $\nu _{T}$ on $\mathcal{B}(\overline{D})$
such that%
\begin{equation}
\mu \left( G\right) =\int_{G}\mathsf{d}\left( \nu _{0}\otimes \nu
_{T}\right) \left( \mathsf{x},\mathsf{y}\right) g(\mathsf{x},T,\mathsf{y})
\label{markovmeasure}
\end{equation}%
for every $G\in \mathcal{B}(\overline{D})\times \mathcal{B}(\overline{D})$,
with $\mu \left( D\times D\right) =1$. It provides therefore a very simple
and practical criterion to decide whether a Bernstein process is Markovian
or not.

\bigskip

It is consequently easy to generate non-Markovian processes out of those
constructed in Theorem 1. One possible way to achieve that and eventually
define the statistical operators and the entropy functional we are
interested in amounts to associating a sequence $\left( p_{\mathsf{m}%
}\right) _{\mathsf{m\in }\mathbb{N}^{d}}$ of probabilities with the pure
point spectrum of (\ref{hamiltonian}), that is, a sequence of numbers
satisfying%
\begin{equation}
p_{\mathsf{m}}\geq 0,\text{ \ \ }\dsum\limits_{\mathsf{m\in }\mathbb{N}%
^{d}}p_{\mathsf{m}}=1,  \label{probabilities}
\end{equation}%
and to consider weighted averages of the form 
\begin{equation}
\overline{\mu }\left( G\right) =\sum_{\mathsf{m\in }\mathbb{N}^{d}}p_{%
\mathsf{m}}\mu _{\mathsf{m}}(G)  \label{mixing}
\end{equation}%
where $\mu _{\mathsf{m}}(G)$ is given by (\ref{probabilitymeasures}). We
note that the preceding series is indeed convergent and defines a genuine
probability measure by virtue of (\ref{normalization}) and (\ref%
{probabilities}). However, in order to generate non-Markovian processes from
(\ref{mixing}) we ought to identify its joint probability density in view of
Remark 2. To this end and aside from having the smoothness of $\varphi _{%
\mathsf{m,0}}>0$ and $\psi _{\mathsf{m,}T}>0$ on $\overline{D}$, the
following additional hypothesis turns out to be sufficient:

\bigskip

(H$_{2}$) We have%
\begin{equation*}
\sup_{\mathsf{m\in }\mathbb{N}^{d}}\sup_{\mathsf{x}\in \overline{D}}\varphi
_{\mathsf{m,0}}\left( \mathsf{x}\right) <+\infty
\end{equation*}%
and%
\begin{equation*}
\sup_{\mathsf{m\in }\mathbb{N}^{d}}\sup_{\mathsf{x}\in \overline{D}}\psi _{%
\mathsf{m,}T}\left( \mathsf{x}\right) <+\infty .
\end{equation*}

\bigskip

This hypothesis indeed clearly implies that%
\begin{equation}
\sum_{\mathsf{m\in }\mathbb{N}^{d}}p_{\mathsf{m}}\varphi _{\mathsf{m,0}%
}\left( \mathsf{x}\right) \psi _{\mathsf{m,}T}\left( \mathsf{y}\right)
<+\infty  \label{partialdensity}
\end{equation}%
for all $\mathsf{x,y\in }\overline{D},$ so that the joint probability
density associated with (\ref{mixing}) may be written as%
\begin{equation}
\overline{\mu }\left( \mathsf{x,y}\right) =g(\mathsf{x},T,\mathsf{y})\sum_{%
\mathsf{m\in }\mathbb{N}^{d}}p_{\mathsf{m}}\varphi _{\mathsf{m,0}}(\mathsf{x)%
}\psi _{\mathsf{m,}T}(\mathsf{y}).  \label{jointdensity}
\end{equation}%
Then the following result holds:

\bigskip

\textbf{Theorem 2. }\textit{Assume that Hypotheses (}H$_{1}$\textit{) and (}H%
$_{2}$\textit{) hold, and for every }$\mathsf{m}\in \mathbb{N}^{d}$\textit{\
let }$Z_{\tau \in \left[ 0,T\right] }^{\mathsf{m}}$\textit{\ be the process
of Theorem 1. Let }$\overline{Z}_{\tau \in \left[ 0,T\right] }$\textit{\ be
the Bernstein process obtained by substituting (\ref{mixing}) into the
formulae of Theorem B.1. Then the following statements are valid:}

\textit{(a) The finite-dimensional distributions of the process }$\overline{Z%
}_{\tau \in \left[ 0,T\right] }$ \textit{are }%
\begin{eqnarray*}
&&\mathbb{P}_{\overline{\mu }}\left( \overline{Z}_{t_{1}}\in F_{1},...,%
\overline{Z}_{t_{n}}\in F_{n}\right) \\
&=&\sum_{\mathsf{m\in }\mathbb{N}^{d}}p_{\mathsf{m}}\mathbb{P}_{\mu _{%
\mathsf{m}}}\left( Z_{t_{1}}^{\mathsf{m}}\in F_{1},...,Z_{t_{n}}^{\mathsf{m}%
}\in F_{n}\right)
\end{eqnarray*}%
\textit{for every }$n\in \mathbb{N}^{+}$\textit{\ and all }$%
F_{1},...,F_{n}\in \mathcal{B}(\overline{D})$\textit{, where }$P_{\mu _{%
\mathsf{m}}}\left( Z_{t_{1}}^{\mathsf{m}}\in F_{1},...,Z_{t_{n}}^{\mathsf{m}%
}\in F_{n}\right) $\textit{\ is given either by (\ref{distributionter}) or (%
\ref{distributionquarto}). In addition, if (\ref{partialdensity}) is not of
the form }$\nu _{0}\otimes \nu _{T}$\textit{\ where }$\nu _{0}$\textit{\ and 
}$\nu _{T}$\textit{\ are as in (\ref{markovmeasure}) then }$\overline{Z}%
_{\tau \in \left[ 0,T\right] }$ \textit{is non-Markovian.}

\textit{(b) We have}%
\begin{equation}
\mathbb{P}_{\overline{\mu }}\left( \overline{Z}_{t}\in F\right) =\sum_{%
\mathsf{m\in }\mathbb{N}^{d}}p_{\mathsf{m}}\mathbb{P}_{\mu _{\mathsf{m}%
}}\left( Z_{t}^{\mathsf{m}}\in F\right)  \label{absoluteprobabilitybis}
\end{equation}%
\textit{for each }$t\in \left[ 0,T\right] $\textit{\ and every }$F\in 
\mathcal{B}(\overline{D})$\textit{, where }$P_{\mu _{\mathsf{m}}}\left(
Z_{t}^{\mathsf{m}}\in F\right) $\textit{\ is given by (\ref%
{absoluteprobability}).}

\textit{(c) We have }%
\begin{equation}
\mathbb{E}_{\overline{\mu }}\left( b(\overline{Z}_{t})\right) =\sum_{\mathsf{%
m\in }\mathbb{N}^{d}}p_{\mathsf{m}}\mathbb{E}_{\mu _{\mathsf{m}}}\left(
b(Z_{t}^{\mathsf{m}})\right)  \label{expectationbis}
\end{equation}%
\textit{for each bounded Borel measurable function }$b:\overline{D}\mapsto 
\mathbb{C}$\textit{\ and every }$t\in \left[ 0,T\right] $\textit{, where }$%
E_{\mu _{\mathsf{m}}}\left( b(Z_{t}^{\mathsf{m}})\right) $\textit{\ is given
by (\ref{expectations}). }

\bigskip

\textbf{Proof.} It follows from Theorem B.1 of the Appendix that a Bernstein
process generated from a statistical mixture of probability measures
coincides with the statistical mixture of the processes generated from those
measures, so that Theorem 2 follows immediately from Theorem 1 and (\ref%
{mixing}). The fact that the process $\overline{Z}_{\tau \in \left[ 0,T%
\right] }$ is non-Markovian when the structural hypothesis regarding (\ref%
{partialdensity}) holds is a direct consequence of Remark 2. \ \ $%
\blacksquare $

\bigskip

\textsc{Remark.} The structural hypothesis we just referred to is necessary
in that it allows one to disregard cases like $\varphi _{\mathsf{m,0}%
}=\varphi _{\mathsf{0}}$ or $\psi _{\mathsf{m,}T}=$ $\psi _{T}$ for every $%
\mathsf{m}$, or the situation where $p_{\mathsf{m}^{\ast }}=1$\textit{\ }for
some\textit{\ }$\mathsf{m}^{\ast }\in \mathbb{N}^{d}$, among others. Indeed,
initial or final data that are identical for each level of the spectrum
still lead to a joint density like that of (\ref{markovmeasure}) with $\nu
_{0}=\varphi _{\mathsf{0}}$ or $\nu _{T}=\psi _{T}$, and hence to a
Markovian process as is the case when $p_{\mathsf{m}^{\ast }}=1$\textit{\ }%
for some\textit{\ }$\mathsf{m}^{\ast }\in \mathbb{N}^{d}$. We shall dwell a
bit more on this further below when we deal with the example in Section 3.

\bigskip We now enquire about the possibility of choosing 
\begin{equation}
\left\{ 
\begin{array}{c}
\varphi _{\mathsf{m,0}}=\mathsf{f}_{\mathsf{m}}, \\ 
\\ 
\psi _{\mathsf{m,}T}=\exp \left[ T\lambda _{\mathsf{m}}\right] \mathsf{f}_{%
\mathsf{m}}%
\end{array}%
\right.  \label{initialfinaldata}
\end{equation}%
as initial-final-data in (\ref{forwardsolution}) and (\ref{backwardsolution}%
), where $\left( \lambda _{\mathsf{m}}\right) _{\mathsf{m}\in \mathbb{N}%
^{d}} $ and $\left( \mathsf{f}_{\mathsf{m}}\right) _{\mathsf{m}\in \mathbb{N}%
^{d}}$ stand for the eigenvalues and eigenfunctions of (\ref{hamiltonian}),
respectively. The difficulty is that the eigenfunctions $\mathsf{f}_{\mathsf{%
m}}$ are not positive in general with the possible exception of $\mathsf{f}%
_{0}$, so that the $\mu _{\mathsf{m}}$ are no longer positive measures with
the possible exception of $\mu _{\mathsf{0}}$. Therefore, we may not
associate a Bernstein process with each level of the spectrum as we did in
Theorem 1. Nevertheless, we proceed by showing that the above averaging
method still allows us to get genuine probability measures in certain cases.
We begin by proving that the $\mu _{\mathsf{m}}$ satisfy the correct
normalization condition under an additional hypothesis:

\bigskip

\textbf{Lemma 1.} \textit{For each }$\mathsf{m}\in \mathbb{N}^{d}$, \textit{%
let us consider measures }$\mu _{\mathsf{m}}$\textit{\ of the form (\ref%
{probabilitymeasures}) where }$\varphi _{\mathsf{m,0}}$\textit{\ and }$\psi
_{\mathsf{m,}T}$\textit{\ are given by (\ref{initialfinaldata}). Then }$\mu
_{\mathsf{m}}$\textit{\ is a signed measure. Moreover, if }%
\begin{equation}
\mathcal{Z(}T):=\sum_{\mathsf{n}\in \mathbb{N}^{d}}\exp \left[ -T\lambda _{%
\mathsf{n}}\right] <+\infty  \label{partitionfunction}
\end{equation}%
\textit{we have}%
\begin{equation}
\mu _{\mathsf{m}}\left( D\times D\right) =1.  \label{normalizationbis}
\end{equation}

\bigskip

\textbf{Proof.} We have just explained why $\mu _{\mathsf{m}}$ is not a
positive measure, so that we need only prove (\ref{normalizationbis}). Since
(\ref{partitionfunction}) holds we have the spectral decomposition%
\begin{equation}
g(\mathsf{x},T,\mathsf{y})=\sum_{\mathsf{n}\in \mathbb{N}^{d}}\exp \left[
-T\lambda _{\mathsf{n}}\right] \mathsf{f}_{\mathsf{n}}(\mathsf{x)f}_{\mathsf{%
n}}(\mathsf{y})  \label{heatkernelbis}
\end{equation}%
as a strongly convergent series in $L^{2}(D\times D)$ for heat kernel (\ref%
{heatkernel}). Therefore, from (\ref{probabilitymeasures}) and (\ref%
{initialfinaldata}) we obtain%
\begin{eqnarray*}
\mu _{\mathsf{m}}(D\times D) &=&\exp \left[ T\lambda _{\mathsf{m}}\right]
\int_{D\times D}\mathsf{dxdyf}_{\mathsf{m}}(\mathsf{x)}g(\mathsf{x},T,%
\mathsf{y})\mathsf{f}_{\mathsf{m}}(\mathsf{y}) \\
&=&\sum_{\mathsf{n}\in \mathbb{N}^{d}}\exp \left[ T(\lambda _{\mathsf{m}%
}-\lambda _{\mathsf{n}})\right] \left( \mathsf{f}_{\mathsf{m}},\mathsf{f}_{%
\mathsf{n}}\right) _{2}^{2}=1
\end{eqnarray*}%
as a consequence of the orthogonality properties of $\left( \mathsf{f}_{%
\mathsf{m}}\right) _{\mathsf{m}\in \mathbb{N}^{d}}$. \ \ $\blacksquare $

\bigskip

Sequences of Gibbs probabilities of the form%
\begin{equation}
p_{\mathsf{m}}=\mathcal{Z}^{-1}(T)\exp \left[ -T\lambda _{\mathsf{m}}\right]
\label{gibbs}
\end{equation}%
will play an important r\^{o}le in the sequel. In fact, with (\ref{gibbs})
the joint probability density of the statistical mixture of the $\mu _{%
\mathsf{m}}$ in Lemma 1 reads%
\begin{eqnarray}
\overline{\mu }(\mathsf{x,y}) &=&g(\mathsf{x},T,\mathsf{y})\sum_{\mathsf{%
m\in }\mathbb{N}^{d}}p_{\mathsf{m}}\exp \left[ T\lambda _{\mathsf{m}}\right] 
\mathsf{f}_{\mathsf{m}}(\mathsf{x)f}_{\mathsf{m}}(\mathsf{y})  \notag \\
&=&\mathcal{Z}^{-1}(T)g(\mathsf{x},T,\mathsf{y})\sum_{\mathsf{m\in }\mathbb{N%
}^{d}}\mathsf{f}_{\mathsf{m}}(\mathsf{x)f}_{\mathsf{m}}(\mathsf{y})  \notag
\\
&=&\mathcal{Z}^{-1}(T)g(\mathsf{x},T,\mathsf{y})\delta (\mathsf{x}-\mathsf{y}%
)  \label{nonmarkovmeasure}
\end{eqnarray}%
as a consequence of the completeness of the basis $\left( \mathsf{f}_{%
\mathsf{m}}\right) _{\mathsf{m}\in \mathbb{N}^{d}}.$ Thus, having (\ref%
{heatkernel}) and (\ref{nonmarkovmeasure}) at our disposal, the latter
obviously not being of the form (\ref{markovmeasure}), and substituting (\ref%
{nonmarkovmeasure}) into Theorem B.1 of Appendix B we obtain:

\bigskip

\textbf{Theorem 3. }\textit{Let us assume that Hypothesis }(H$_{1}$)\textit{%
\ holds, and let }$\bar{Z}_{\tau \in \left[ 0,T\right] }$ \textit{be the
Bernstein generated by (\ref{nonmarkovmeasure}). Then the following
statements are valid:}

\textit{(a) The process }$\bar{Z}_{\tau \in \left[ 0,T\right] }$ \textit{is
stationary, non-Markovian and for every }$n\in \mathbb{N}^{+}$ \textit{with} 
$n\geq 2$ \textit{its finite-dimensional distributions are}%
\begin{eqnarray}
&&\mathbb{P}_{\bar{\mu}}\left( \bar{Z}_{t_{1}}\in F_{1},...,\bar{Z}%
_{t_{n}}\in F_{n}\right)  \notag \\
&=&\mathcal{Z}^{-1}(T)\int_{F_{1}}\mathsf{dx}_{1}...\int_{F_{n}}\mathsf{dx}%
_{n}  \notag \\
&&\times \dprod\limits_{k=2}^{n}g\left( \mathsf{x}_{k},t_{k}-t_{k-1},\mathsf{%
x}_{k-1}\right) \times g\left( \mathsf{x}_{1},T-(t_{n}-t_{1}),\mathsf{x}%
_{n}\right)  \label{distributionquinto}
\end{eqnarray}%
\textit{for all }$F_{1},...,F_{n}\in \mathcal{B}(\overline{D})$\textit{\ and
all }$0<t_{1}<...<t_{n}<T$\textit{. }

\textit{(b) We have}%
\begin{equation}
\mathbb{P}_{\bar{\mu}}\left( \bar{Z}_{t}\in F\right) =\mathcal{Z}%
^{-1}(T)\int_{F}\mathsf{dx}g\left( \mathsf{x},T,\mathsf{x}\right)
\label{probability}
\end{equation}%
\textit{for each }$F\in \mathcal{B}(\overline{D})$\textit{\ and every} $t\in %
\left[ 0,T\right] $.

\bigskip \textit{(c) We have}%
\begin{equation}
\mathbb{E}_{_{\bar{\mu}}}\left( b(\bar{Z}_{t})\right) =\mathcal{Z}%
^{-1}(T)\int_{D}\mathsf{dx}b(\mathsf{x})g\left( \mathsf{x},T,\mathsf{x}%
\right)  \label{expectationter}
\end{equation}%
\textit{for each bounded Borel measurable function }$b:\overline{D}\mathbb{%
\mapsto C}$ \textit{and every }$t\in \left[ 0,T\right] $\textit{.}

\bigskip

\textsc{Remark.} The fact that the process of the preceding result is
stationary is tied up with the structure of the finite-dimensional
distributions (\ref{distributionquinto}), which differs from those in
Theorems 1 and 2. Indeed, for any $\tau >0$ sufficiently small such that $%
0<t_{1}+\tau <...<t_{n}+\tau <T$ we have%
\begin{equation*}
\mathbb{P}_{\bar{\mu}}\left( \bar{Z}_{t_{1}+\tau }\in F_{1},...,\bar{Z}%
_{t_{n}+\tau }\in F_{n}\right) =\mathbb{P}_{\bar{\mu}}\left( \bar{Z}%
_{t_{1}}\in F_{1},...,\bar{Z}_{t_{n}}\in F_{n}\right) ,
\end{equation*}%
as well as the time independence of (\ref{probability}) and (\ref%
{expectationter}). Furthermore we also note that since $\mathbb{P}_{\bar{\mu}%
}\left( \bar{Z}_{t}\in D\right) =1$, Relation (\ref{probability}) provides
yet another expression for (\ref{partitionfunction}), namely%
\begin{equation*}
\mathcal{Z}(T)=\int_{D}\mathsf{dx}g\left( \mathsf{x},T,\mathsf{x}\right)
\end{equation*}%
which, of course, also follows from (\ref{heatkernelbis}) and the fact that $%
\left\Vert \mathsf{f}_{\mathsf{m}}\right\Vert _{2}=1$ for every $\mathsf{m}%
\in \mathbb{N}^{d}$.

\bigskip

The preceding results thus reveal the possibility of having at least two
types of Bernstein processes, namely, on the one hand Markovian processes
associated with each level of the spectrum of (\ref{hamiltonian}), and on
the other hand typically non-Markovian processes obtained by averaging $%
Z_{\tau \in \left[ 0,T\right] }^{\mathsf{m}}$ over the whole spectrum for a
given sequence $\left( p_{\mathsf{m}}\right) _{\mathsf{m\in }\mathbb{N}^{d}}$%
, or by averaging signed measures. In order to better characterize those
processes by means of entropy considerations, we now proceed by introducing
a statistical operator and an entropy functional by analogy with Quantum
Statistical Mechanics. We define%
\begin{equation}
\mathcal{R}f:=\sum_{\mathsf{m\in }\mathbb{N}^{d}}p_{\mathsf{m}}\left( f,%
\mathsf{f}_{\mathsf{m}}\right) _{2}\mathsf{f}_{\mathsf{m}}
\label{statisticaloperator}
\end{equation}%
for each $f\in L^{2}(D)$. The following result is elementary, so that we
only sketch the proof of the trace-class property which will be discussed in
a more general context in Appendix A:

\bigskip

\textbf{Proposition 1.} \textit{Let us assume that Hypothesis (}H$_{1}$%
\textit{) holds. Then the following statements are valid:}

\textit{(a) Expression (\ref{statisticaloperator}) defines a self-adjoint,
positive trace-class operator in }$L^{2}\left( D\right) $\textit{\ such
that\ the inequalities}%
\begin{equation*}
0\leq \mathcal{R}^{2}\leq \mathcal{R}\leq \mathbb{I}
\end{equation*}%
\textit{hold} \textit{in the sense of quadratic forms, where }$\mathbb{I}$%
\textit{\ stands for the identity in }$L^{2}\left( D\right) $\textit{. More
specifically we have}%
\begin{equation}
\func{Tr}\mathcal{R}=\sum_{\mathsf{m\in }\mathbb{N}^{d}}p_{\mathsf{m}}=1
\label{trace}
\end{equation}%
\textit{and}%
\begin{equation}
\func{Tr}\mathcal{R}^{2}=\sum_{\mathsf{m\in }\mathbb{N}^{d}}p_{\mathsf{m}%
}^{2}\leq 1.  \label{tracebis}
\end{equation}%
\textit{In particular we have}%
\begin{equation*}
\func{Tr}\mathcal{R}^{2}=1
\end{equation*}%
\textit{\ if, and only if, }$p_{\mathsf{m}^{\ast }}=1$\textit{\ for some }$%
\mathsf{m}^{\ast }\in \mathbb{N}^{d}$\textit{\ and thus }$p_{\mathsf{m}}=0$%
\textit{\ for every }$\mathsf{m\neq m}^{\ast }$\textit{.}

\textit{(b) The eigenvalue equation}%
\begin{equation}
\mathcal{R}\mathsf{f}_{\mathsf{m}}=p_{\mathsf{m}}\mathsf{f}_{m}
\label{eigenvalueequation}
\end{equation}%
\textit{holds for each }$\mathsf{m}\in \mathbb{N}^{d}$\textit{\ and the
spectrum of }$\mathcal{R}$ \textit{is either pure point with }$\sigma (%
\mathcal{R})=(p_{\mathsf{m}})_{\mathsf{m\in }\mathbb{N}^{d}\text{ }}$\textit{%
if }$p_{\mathsf{m}}=0$ \textit{for at least one} $\mathsf{m}$, \textit{or }$%
\sigma (\mathcal{R})=(p_{\mathsf{m}})_{\mathsf{m\in }\mathbb{N}^{d}\text{ }%
}\cup \left\{ 0\right\} $ \textit{if} $0<p_{\mathsf{m}}<1$ \textit{for every}
$\mathsf{m}$,\textit{\ in which case zero is not an eigenvalue.}

\textit{(c) If }$B$\textit{\ is a linear bounded self-adjoint operator on }$%
L^{2}(D)$ \textit{we have}%
\begin{equation}
\func{Tr}\left( \mathcal{R}B\right) =\sum_{\mathsf{m\in }\mathbb{N}^{d}}p_{%
\mathsf{m}}\left( B\mathsf{f}_{\mathsf{m}},\mathsf{f}_{\mathsf{m}}\right)
_{2}.  \label{stataeveragebis}
\end{equation}%
\textit{In particular, if }$B$\textit{\ is the multiplication operator given
by} $Bf=bf$ \textit{where }$b\in L^{\infty }(D)$\textit{\ is real-valued, we
have}%
\begin{equation}
\func{Tr}\left( \mathcal{R}B\right) \mathcal{=}\sum_{\mathsf{m\in }\mathbb{N}%
^{d}}p_{\mathsf{m}}\int_{D}d\mathsf{x}b(\mathsf{x})\left\vert \mathsf{f}_{%
\mathsf{m}}(\mathsf{x})\right\vert ^{2}.  \label{stataverage}
\end{equation}

\bigskip

\textbf{Proof. }Owing to the properties of $p_{\mathsf{m}}$ and $\mathsf{f}_{%
\mathsf{m}}$ it is immediate that (\ref{statisticaloperator}) defines a
linear bounded operator in $L^{2}(D)$. Now let $\left( \mathsf{h}_{\mathsf{n}%
}\right) _{\mathsf{n\in }\mathbb{N}^{d}}$ be an arbitrary orthonormal basis
in $L^{2}(D)$. In order to prove that $\mathcal{R}$ is trace-class, it is
then necessary and sufficient to show that 
\begin{equation*}
\sum_{\mathsf{n\in }\mathbb{N}^{d}}\left( \mathcal{R}\mathsf{h}_{\mathsf{n}},%
\mathsf{h}_{\mathsf{n}}\right) _{2}<+\infty
\end{equation*}%
(see, e.g., Theorem 8.1 in Chapter III of \cite{gohbergkrein}). To this end
let us introduce the function%
\begin{equation}
a(\mathsf{m},\mathsf{n}):=p_{\mathsf{m}}\left( \mathsf{h}_{\mathsf{n}},%
\mathsf{f}_{\mathsf{m}}\right) _{2}\left( \mathsf{f}_{\mathsf{m}},\mathsf{h}%
_{\mathsf{n}}\right) _{2}  \label{auxiliaryfunction}
\end{equation}%
so that%
\begin{equation}
\sum_{\mathsf{m\in }\mathbb{N}^{d}}a(\mathsf{m},\mathsf{n})=\left( \mathcal{R%
}\mathsf{h}_{\mathsf{n}},\mathsf{h}_{\mathsf{n}}\right) _{2}
\label{expression1}
\end{equation}%
for every fixed $\mathsf{n}$. Moreover, for any fixed $\mathsf{m}$ we have%
\begin{equation}
\sum_{\mathsf{n\in }\mathbb{N}^{d}}a(\mathsf{m},\mathsf{n})=p_{\mathsf{m}}
\label{expression2}
\end{equation}%
since $\left\Vert \mathsf{f}_{\mathsf{m}}\right\Vert _{2}=1$. Furthermore,
the preceding series converges absolutely since from (\ref{auxiliaryfunction}%
) we have for any choice of positive integers $N_{1},...,N_{d}$ the estimate%
\begin{eqnarray*}
&&\sum_{\mathsf{n:0\leq n}_{j}\leq N_{j}}\left\vert a(\mathsf{m},\mathsf{n}%
)\right\vert \\
&\leq &p_{\mathsf{m}}\left( \sum_{\mathsf{n\in }\mathbb{N}^{d}}\left\vert
\left( \mathsf{h}_{\mathsf{n}},\mathsf{f}_{\mathsf{m}}\right)
_{2}\right\vert ^{2}\right) ^{\frac{1}{2}}\left( \sum_{\mathsf{n\in }\mathbb{%
N}^{d}}\left\vert \left( \mathsf{f}_{\mathsf{m}},\mathsf{h}_{\mathsf{n}%
}\right) _{2}\right\vert ^{2}\right) ^{\frac{1}{2}}=p_{\mathsf{m}}
\end{eqnarray*}%
for any fixed $\mathsf{m}$. Consequently we have 
\begin{equation*}
\sum_{\mathsf{m\in }\mathbb{N}^{d}}\sum_{\mathsf{n\in }\mathbb{N}%
^{d}}\left\vert a(\mathsf{m},\mathsf{n})\right\vert \leq \sum_{\mathsf{m\in }%
\mathbb{N}^{d}}p_{\mathsf{m}}=1,
\end{equation*}%
from which we infer according to well-known criterias that the associated
iterated series are equal, that is,%
\begin{equation*}
\sum_{\mathsf{n\in }\mathbb{N}^{d}}\sum_{\mathsf{m\in }\mathbb{N}^{d}}a(%
\mathsf{m},\mathsf{n})=\sum_{\mathsf{m\in }\mathbb{N}^{d}}\sum_{\mathsf{n\in 
}\mathbb{N}^{d}}a(\mathsf{m},\mathsf{n}).
\end{equation*}%
Equivalently, this means that%
\begin{equation*}
\func{Tr}\mathcal{R}:=\sum_{\mathsf{n\in }\mathbb{N}^{d}}\left( \mathcal{R}%
\mathsf{h}_{\mathsf{n}},\mathsf{h}_{\mathsf{n}}\right) _{2}=\sum_{\mathsf{%
m\in }\mathbb{N}^{d}}p_{\mathsf{m}}=1
\end{equation*}%
according to (\ref{expression1}) and (\ref{expression2}), which proves (\ref%
{trace}). The proof of (\ref{tracebis}) is similar with%
\begin{equation*}
\mathcal{R}^{2}f=\sum_{\mathsf{m\in }\mathbb{N}^{d}}p_{\mathsf{m}}^{2}\left(
f,\mathsf{f}_{\mathsf{m}}\right) _{2}\mathsf{f}_{\mathsf{m}}.
\end{equation*}%
The remaining statements are immediate from elementary arguments. \ \ $%
\blacksquare $

\bigskip

\textsc{Remark.} Regarding expression (\ref{stataverage}) we note that when
the $p_{\mathsf{m}}$ are given by (\ref{gibbs}) we have 
\begin{equation}
\func{Tr}\left( \mathcal{R}B\right) =\mathbb{E}_{_{\bar{\mu}}}\left( b(\bar{Z%
}_{t})\right)  \label{stataverageter}
\end{equation}%
for every $t\in \left[ 0,T\right] $, where the right-hand side is given by (%
\ref{expectationter}). This is an immediate consequence of (\ref%
{heatkernelbis}), so that the statistical average (\ref{stataverage})
calculated by means of Gibbs probabilities coincides with the expectation of
some function of the process of Theorem 3. This is of course only possible
because that process is stationary, the right-hand side of (\ref%
{stataverageter}) then being time-independent as the left-hand side is. It
is therefore reasonable to ask whether relations such as (\ref%
{stataverageter}) may exist in more general cases, for instance for the
averaged processes of Theorem 2 which are in general non-stationary. This is
indeed possible as we shall show in the appendix, provided we have at our
disposal a class of time-dependent statistical operators which generalize (%
\ref{statisticaloperator}).

\bigskip

By analogy with Quantum Statistical Mechanics from which we also borrow the
terminology (see, e.g., Section 3 in Chapter V of \cite{vonneumann}),
Proposition 1 allows us to establish a preliminary classification of the
Bernstein processes constructed above, according to the following:

\bigskip

\textbf{Definition 2. }\textit{For a given sequence }$\left( p_{\mathsf{m}%
}\right) _{\mathsf{m\in }\mathbb{N}^{d}}$\textit{\ let }$\overline{Z}_{\tau
\in \left[ 0,T\right] }$\textit{\ be the Bernstein process of Theorem 2, and
let }$\mathcal{R}$ \textit{be the statistical operator given by (\ref%
{statisticaloperator}).}

\textit{(a) If }$\func{Tr}\mathcal{R}^{2}=1$\textit{\ we say that }$%
\overline{Z}_{\tau \in \left[ 0,T\right] }$\textit{\ is a pure process.}

\textit{(b) If }$\func{Tr}\mathcal{R}^{2}<1$\textit{\ we say that }$%
\overline{Z}_{\tau \in \left[ 0,T\right] }$\textit{\ is a mixed process.}

\bigskip

We note that in the first case we necessarily have $\overline{Z}_{\tau \in %
\left[ 0,T\right] }=Z_{\tau \in \left[ 0,T\right] }^{\mathsf{m}^{\ast }}$
for some $\mathsf{m}^{\ast }\in \mathbb{N}^{d}$ according to the second part
of (a) in Proposition 1, so that $\overline{Z}_{\tau \in \left[ 0,T\right] }$
reduces to a Markovian process according to Theorem 1 or the remark
following the proof of Theorem 2. On the other hand, an important example
which illustrates the second case is that of Gibbs probability measures (\ref%
{gibbs}).

\bigskip

We now introduce the entropy functional

\begin{equation}
\mathsf{S:}=\sum_{\mathsf{m\in }\mathbb{N}^{d}}p_{\mathsf{m}}\ln \left( 
\frac{1}{p_{\mathsf{m}}}\right)  \label{entropyfunctional}
\end{equation}%
where we define $x\ln \left( \frac{1}{x}\right) $ to be zero at $x=0$ so
that $\mathsf{S=0}$ if, and only if, $p_{\mathsf{m}}=0$ or $p_{\mathsf{m}}=1$
for every $\mathsf{m}$, the latter value being associated with pure
processes according to Definition 2. It is plain that we may have $\mathsf{%
S=+\infty }$ despite the normalization (\ref{probabilities}), a case in
point being that of the Gibbs probabilities (\ref{gibbs}). Indeed, the
substitution of (\ref{gibbs}) into (\ref{entropyfunctional}) shows that $%
\mathsf{S<+\infty }$ if, and only if, the additional condition%
\begin{equation*}
\sum_{\mathsf{m}\in \mathbb{N}^{d}}\exp \left[ -T\lambda _{\mathsf{m}}\right]
\lambda _{\mathsf{m}}<+\infty
\end{equation*}%
holds. From now on we shall therefore assume that the $p_{\mathsf{m}}$ are
chosen in such a way that $0<p_{\mathsf{m}}<1$ with%
\begin{equation}
\sum_{\mathsf{m\in }\mathbb{N}^{d}}p_{\mathsf{m}}\ln \left( \frac{1}{p_{%
\mathsf{m}}}\right) <+\infty .  \label{finitude}
\end{equation}%
The following result is then our desired optimization statement for (\ref%
{entropyfunctional}). We note that we only consider there probabilities
which assign an \textit{a priori} prescribed value to the average of the
spectrum of (\ref{hamiltonian}):

\bigskip

\textbf{Theorem 4.}\textit{\ Let us consider the set of all sequences }$%
\left( p_{\mathsf{m}}\right) _{\mathsf{m\in }\mathbb{N}^{d}}$ \textit{%
satisfying }$0<p_{\mathsf{m}}<1$\textit{\ for every }$\mathsf{m}$,\textit{\
along with}%
\begin{equation}
\sum_{\mathsf{m\in }\mathbb{N}^{d}}p_{\mathsf{m}}=1  \label{normalizationter}
\end{equation}%
\textit{\ and (\ref{finitude})}. \textit{Moreover}, \textit{let }$\lambda
\in \mathbb{R}$\textsf{\ }\textit{be given and let us assume that}%
\begin{equation}
\sum_{\mathsf{m\in }\mathbb{N}^{d}}p_{\mathsf{m}}\lambda _{\mathsf{m}%
}=\lambda .  \label{averagespectrum}
\end{equation}%
\textit{Then the following statements are valid:}

\textit{(a) There exists a finite constant }$\beta (\lambda )>0$ \textit{%
such that}%
\begin{equation}
\mathcal{Z}(\beta ):=\sum_{\mathsf{m}\in \mathbb{N}^{d}}\exp \left[ -\beta
\lambda _{\mathsf{m}}\right] <+\infty  \label{partitionfunctionbis}
\end{equation}%
\textit{and}%
\begin{equation}
\mathcal{Z}^{-1}(\beta )\sum_{\mathsf{m}\in \mathbb{N}^{d}}\exp \left[
-\beta \lambda _{\mathsf{m}}\right] \lambda _{\mathsf{m}}<+\infty
\label{averagespectrumbis}
\end{equation}%
\textit{for every} $\beta \in \left[ \beta (\lambda ),+\infty \right) $.

\textit{(b) Among all the mixed processes obtained from sequences of the
above type by the method of Theorem 2, the process of maximal entropy is
that generated from probabilities given by}%
\begin{equation}
p_{\mathsf{m,Gibbs}}=\mathcal{Z}^{-1}(\beta (\lambda ))\exp \left[ -\beta
(\lambda )\lambda _{\mathsf{m}}\right]  \label{gibbsbis}
\end{equation}%
\textit{for every} $\mathsf{m}\in \mathbb{N}^{d}$. \textit{Moreover we have}%
\begin{equation}
\mathcal{Z}^{-1}(\beta (\lambda ))\sum_{\mathsf{m}\in \mathbb{N}^{d}}\exp %
\left[ -\beta (\lambda )\lambda _{\mathsf{m}}\right] \lambda _{\mathsf{m}%
}=\lambda .  \label{averagespectrumter}
\end{equation}

\textit{(c) If we assume in addition that }$\sum_{\mathsf{m}\in \mathbb{N}%
^{d}}\exp \left[ -\beta \lambda _{\mathsf{m}}\right] \lambda _{\mathsf{m}%
}<+\infty $ \textit{for every} $\beta \in \left( 0,\beta (\lambda )\right) $,%
\textit{\ then} $\beta \mapsto \mathcal{Z}(\beta )$ \textit{is\
differentiable at }$\beta =\beta (\lambda )$\textit{\ and we have}%
\begin{equation}
\mathsf{S}_{\max }\left( \lambda \right) =\ln \mathcal{Z}(\beta (\lambda
))-\beta (\lambda )\frac{d}{d\beta }\ln \mathcal{Z}(\beta )_{\left\vert
\beta =\beta (\lambda )\right. }  \label{maxentropy}
\end{equation}%
\textit{for the maximal entropy of part (b). }

\bigskip

\textbf{Proof.} Since $\lambda _{\mathsf{m}}\rightarrow +\infty $ as $%
\left\vert \mathsf{m}\right\vert \rightarrow +\infty $ and since (\ref%
{normalizationter}) holds, there exist $\mathsf{n,n}^{\prime }\in $ $\mathbb{%
N}^{d}$ with $\mathsf{n\neq n}^{\prime }$ such that $\lambda _{\mathsf{n}%
}\neq \lambda _{\mathsf{n}^{\prime }}$ and $p_{\mathsf{n}}\neq p_{\mathsf{n}%
^{\prime }}$. We then consider the inhomogeneous system%
\begin{eqnarray}
\alpha +\beta \lambda _{\mathsf{n}} &=&-\left( \ln p_{\mathsf{n}}+1\right) ,
\label{equation1} \\
\alpha +\beta \lambda _{\mathsf{n}^{\prime }} &=&-\left( \ln p_{\mathsf{n}%
^{\prime }}+1\right)  \label{equation2}
\end{eqnarray}%
in the two unknowns $\alpha $ and $\beta $, whose unique solution pair reads%
\begin{eqnarray}
\alpha &=&\left( \lambda _{\mathsf{n}^{\prime }}-\lambda _{\mathsf{n}%
}\right) ^{-1}\left( \lambda _{\mathsf{n}}\left( \ln p_{\mathsf{n}^{\prime
}}+1\right) -\lambda _{\mathsf{n}^{\prime }}\left( \ln p_{\mathsf{n}%
}+1\right) \right) ,  \label{expre1} \\
\beta &=&\left( \lambda _{\mathsf{n}^{\prime }}-\lambda _{\mathsf{n}}\right)
^{-1}\ln \frac{p_{\mathsf{n}}}{p_{\mathsf{n}^{\prime }}}.  \label{expre2}
\end{eqnarray}%
Furthermore, let us write (\ref{normalizationter}) and (\ref{averagespectrum}%
) as%
\begin{eqnarray*}
p_{\mathsf{n}}+p_{\mathsf{n}^{\prime }} &=&1-\sum_{\mathsf{m\in }\mathbb{N}%
^{d},\text{ }\mathsf{m}\neq \mathsf{n,}\text{ }\mathsf{n}^{\prime }}p_{%
\mathsf{m}}, \\
p_{\mathsf{n}}\lambda _{\mathsf{n}}+p_{\mathsf{n}^{\prime }}\lambda _{%
\mathsf{n}^{\prime }} &=&\lambda -\sum_{\mathsf{m\in }\mathbb{N}^{d},\text{ }%
\mathsf{m}\neq \mathsf{n,}\text{ }\mathsf{n}^{\prime }}p_{\mathsf{m}}\lambda
_{\mathsf{m}},
\end{eqnarray*}%
respectively, which gives%
\begin{equation*}
p_{\mathsf{n}}=\left( \lambda _{\mathsf{n}^{\prime }}-\lambda _{\mathsf{n}%
}\right) ^{-1}\left( \lambda _{\mathsf{n}^{\prime }}\left( 1-\sum_{\mathsf{%
m\in }\mathbb{N}^{d},\text{ }\mathsf{m}\neq \mathsf{n,}\text{ }\mathsf{n}%
^{\prime }}p_{\mathsf{m}}\right) -\lambda +\sum_{\mathsf{m\in }\mathbb{N}%
^{d},\text{ }\mathsf{m}\neq \mathsf{n,}\text{ }\mathsf{n}^{\prime }}p_{%
\mathsf{m}}\lambda _{\mathsf{m}}\right)
\end{equation*}%
and 
\begin{equation*}
p_{\mathsf{n}^{\prime }}=\left( \lambda _{\mathsf{n}^{\prime }}-\lambda _{%
\mathsf{n}}\right) ^{-1}\left( \lambda -\sum_{\mathsf{m\in }\mathbb{N}^{d},%
\text{ }\mathsf{m}\neq \mathsf{n,}\text{ }\mathsf{n}^{\prime }}p_{\mathsf{m}%
}\lambda _{\mathsf{m}}-\lambda _{\mathsf{n}}\left( 1-\sum_{\mathsf{m\in }%
\mathbb{N}^{d},\text{ }\mathsf{m}\neq \mathsf{n,}\text{ }\mathsf{n}^{\prime
}}p_{\mathsf{m}}\right) \right) .
\end{equation*}%
The substitution of these expressions into (\ref{expre1}) and (\ref{expre2})
shows that $\alpha =\alpha \left( \lambda \right) $ and $\beta =\beta \left(
\lambda \right) $ depend on $\lambda $, and that%
\begin{equation}
\ln p_{\mathsf{m}}+1=-\alpha \left( \lambda \right) -\beta \left( \lambda
\right) \lambda _{\mathsf{m}}  \label{expression}
\end{equation}%
for $\mathsf{m=n}$ and $\mathsf{m=n}^{\prime }$ according to (\ref{equation1}%
) and (\ref{equation2}). Now for every $\mathsf{j\in }\mathbb{N}^{d}$ with $%
\mathsf{j}\neq \mathsf{n,}$ $\mathsf{n}^{\prime }$ we have%
\begin{eqnarray}
\frac{\partial p_{\mathsf{n}}}{\partial p_{\mathsf{j}}} &=&\left( \lambda _{%
\mathsf{n}^{\prime }}-\lambda _{\mathsf{n}}\right) ^{-1}\left( \lambda _{%
\mathsf{j}}-\lambda _{\mathsf{n}^{\prime }}\right) ,  \label{derivative1} \\
\frac{\partial p_{\mathsf{n}^{\prime }}}{\partial p_{\mathsf{j}}} &=&\left(
\lambda _{\mathsf{n}^{\prime }}-\lambda _{\mathsf{n}}\right) ^{-1}\left(
\lambda _{\mathsf{n}}-\lambda _{\mathsf{j}}\right) .  \label{derivative2}
\end{eqnarray}%
Furthermore, let us define%
\begin{equation}
\overline{\mathsf{S}}:=p_{\mathsf{n}}\ln p_{\mathsf{n}}+p_{\mathsf{n}%
^{\prime }}\ln p_{\mathsf{n}^{\prime }}+\sum_{\mathsf{m\in }\mathbb{N}^{d},%
\text{ }\mathsf{m}\neq \mathsf{n,}\text{ }\mathsf{n}^{\prime }}p_{\mathsf{m}%
}\ln p_{\mathsf{m}}  \label{entropyfunctionalbis}
\end{equation}%
where $p_{\mathsf{n}}$ and $p_{\mathsf{n}^{\prime }}$ are given by the above
expressions. From (\ref{derivative1}) and (\ref{derivative2}) followed by
the use of (\ref{expre1}) and (\ref{expre2}) we get%
\begin{eqnarray}
\frac{\partial \overline{\mathsf{S}}}{\partial p_{\mathsf{j}}} &=&\left(
\lambda _{\mathsf{n}^{\prime }}-\lambda _{\mathsf{n}}\right) ^{-1}\left(
\left( \lambda _{\mathsf{j}}-\lambda _{\mathsf{n}^{\prime }}\right) \left(
\ln p_{\mathsf{n}}+1\right) +\left( \lambda _{\mathsf{n}}-\lambda _{\mathsf{j%
}}\right) \left( \ln p_{\mathsf{n}^{\prime }}+1\right) \right) +\ln p_{%
\mathsf{j}}+1  \notag \\
&=&\alpha \left( \lambda \right) +\beta \left( \lambda \right) \lambda _{%
\mathsf{j}}+\ln p_{\mathsf{j}}+1,  \label{gradient}
\end{eqnarray}%
so that $\frac{\partial \overline{\mathsf{S}}}{\partial p_{\mathsf{j}}}=0$
if, and only if, 
\begin{equation*}
\ln p_{\mathsf{j}}+1=-\alpha \left( \lambda \right) -\beta \left( \lambda
\right) \lambda _{\mathsf{j}}.
\end{equation*}%
We now combine this with (\ref{expression}) to conclude that for every
choice of probabilities which satisfy the hypotheses of the theorem and
which annihilate (\ref{gradient}) we have%
\begin{equation*}
\ln p_{\mathsf{m}}+1=-\alpha \left( \lambda \right) -\beta \left( \lambda
\right) \lambda _{\mathsf{m}}
\end{equation*}%
for every $\mathsf{m}\in \mathbb{N}^{d}$, that is,%
\begin{equation*}
p_{\mathsf{m}}=\exp \left[ -\left( 1+\alpha \left( \lambda \right) \right) %
\right] \exp \left[ -\beta \left( \lambda \right) \lambda _{\mathsf{m}}%
\right] .
\end{equation*}%
Consequently, since $\lambda _{\mathsf{m}}\rightarrow +\infty $ as $%
\left\vert \mathsf{m}\right\vert \rightarrow +\infty $ we obtain (\ref%
{partitionfunctionbis}) and (\ref{averagespectrumbis}) from (\ref%
{normalizationter}) and (\ref{averagespectrum}), respectively, where $\beta
\left( \lambda \right) >0$ and $\beta \in \left[ \beta (\lambda ),+\infty
\right) $. This proves Statement (a) and gives (\ref{gibbsbis}) for every $%
\mathsf{m}\in \mathbb{N}^{d}$ along with (\ref{averagespectrumter}). The
fact that (\ref{entropyfunctional}) is indeed maximal at (\ref{gibbsbis})
then follows from an adaptation of well-known considerations (see, e.g.,
Section 8 in Chapter 4 of \cite{brillouin}).

The differentiability of $\beta \mapsto \mathcal{Z}(\beta )$ at $\beta
=\beta (\lambda )$ under the stated conditions as well as (\ref{maxentropy})
are consequences of elementary arguments and of the direct substitution of (%
\ref{gibbsbis}) into (\ref{entropyfunctional}). $\ \ \blacksquare $

\bigskip

\textsc{Remarks.} (1) Relation (\ref{averagespectrum}) has to be viewed as a
further restriction on the class of admissible probabilities, which of
course must be such that $p_{\mathsf{m}}\lambda _{\mathsf{m}}\rightarrow 0$
as $\left\vert \mathsf{m}\right\vert \rightarrow +\infty $. Furthermore, the
choice of the preassigned value $\lambda $ must be consistent with the
nature of those eigenvalues. Thus if for instance $\lambda _{\mathsf{m}}\geq
0$ for every $\mathsf{m}$, one must then impose $\lambda \in \mathbb{R}^{+}$
for (\ref{averagespectrum}) to make sense. We present an example of that
kind in Section 3, where we also have $\sum_{\mathsf{m}\in \mathbb{N}%
^{d}}\exp \left[ -\beta \lambda _{\mathsf{m}}\right] \lambda _{\mathsf{m}%
}<+\infty $ for every $\beta \in \left( 0,+\infty \right) $.

(2) Whereas the preceding considerations describe a situation where (\ref%
{entropyfunctional}) does not depend explicitly on time, there are many
time-dependent entropy functionals which we may associate with Bernstein
processes, for instance%
\begin{equation}
\mathsf{S}_{\mathsf{m}}(t):=\int_{D\times D}d\mathsf{x}d\mathsf{y}w_{\mathsf{%
m}}(\mathsf{x},t;\mathsf{y},0)\ln \left( \frac{1}{w_{\mathsf{m}}(\mathsf{x}%
,t;\mathsf{y},0)}\right)  \label{entropyfunctionalter}
\end{equation}%
in case of the Markovian processes of Theorem 1, where $w_{\mathsf{m}}$
given by (\ref{markovdensitybis}) satisfies a specific Kolmogorov or
Fokker-Planck equation. We defer the derivation of such equations, the
analyses of the related entropy functionals such as (\ref%
{entropyfunctionalter}) and their consequences to a separate publication.

\bigskip

We devote the next section to illustrating some of the above results.

\section{A hierarchy of Bernstein processes in a two-dimensional disk}

We consider here forward-backward problems of the form (\ref{forwardproblem}%
)-(\ref{backwardproblem}) with $V=0$ identically, defined in the open
two-dimensional disk of radius one centered at the origin, so that
Hypothesis (H$_{1}$) trivially holds. We limit ourselves to an illustration
of a few properties listed in the preceding section regarding Bernstein
processes generated by certain radially symmetric solutions to such
problems. Thus, we first switch to polar coordinates and start out with the
hierarchy%
\begin{align}
\partial _{t}u(r,t)& =\frac{1}{2}\left( \frac{\partial ^{2}}{\partial r^{2}}+%
\frac{1}{r}\frac{\partial }{\partial r}\right) u(r,t),\text{ \ }(r,t)\in
\left( 0,1\right] \times \left( 0,T\right] ,  \notag \\
u(r,0)& =\varphi _{0,\mathsf{m}}(r),\text{ \ \ }r\in \left[ 0,1\right] ,
\label{forwardproblembis} \\
\partial _{r}u(1,t)& =0,\text{ \ \ \ }t\in \left[ 0,T\right]   \notag
\end{align}%
and%
\begin{align}
-\partial _{t}v(r,t)& =\frac{1}{2}\left( \frac{\partial ^{2}}{\partial r^{2}}%
+\frac{1}{r}\frac{\partial }{\partial r}\right) v(r,t),\text{ \ }(r,t)\in
\left( 0,1\right] \times \left[ 0,T\right) ,  \notag \\
v(r,T)& =\psi _{T,\mathsf{m}}(r),\text{ \ \ }r\in \left[ 0,1\right] ,
\label{backwardproblembis} \\
\partial _{r}v(1,t)& =0,\text{ \ \ }t\in \left[ 0,T\right] .  \notag
\end{align}%
In this case the index $\mathsf{m}\in \mathbb{N}$ labels the discrete
spectrum of the radial part of Neumann's Laplacian $-\frac{1}{2}\Delta _{%
\mathsf{x}}$ on the disk, which consists exclusively of eigenvalues $\lambda
_{\mathsf{m}}\geq 0$ determined by the condition%
\begin{equation}
J_{1}\left( \sqrt{2\lambda _{\mathsf{m}}}\right) =0  \label{neumann}
\end{equation}%
where $J_{1}$ stands for the Bessel function of the first kind of order one.
For convenience we order these eigenvalues as%
\begin{equation}
0=\lambda _{\mathsf{0}}<\lambda _{1}<\lambda _{\mathsf{2}}<.....\text{ ,%
\label{ordering}}
\end{equation}%
and recall that there exists a finite constant $c>0$ such that%
\begin{equation}
c\left( \mathsf{m}-1\right) ^{2}<\lambda _{\mathsf{m}}<c\left( \mathsf{m}%
+1\right) ^{2}  \label{growth}
\end{equation}%
for every $\mathsf{m}\in \mathbb{N}^{+}$. Moreover, the corresponding
orthonormal basis $\left( \mathsf{f}_{\mathsf{m}}\right) _{\mathsf{m}\in 
\mathbb{N}}$ of eigenfunctions in the space of all complex-valued,
square-integrable functions with respect to the measure $rdr$ on $\left(
0,1\right) $ is given by%
\begin{equation}
\mathsf{f}_{\mathsf{m}}\left( r\right) =\frac{\sqrt{2}}{\left\vert
J_{0}\left( \sqrt{2\lambda _{\mathsf{m}}}\right) \right\vert }J_{0}\left( 
\sqrt{2\lambda _{\mathsf{m}}}r\right)   \label{eigenfunctions}
\end{equation}%
where $J_{0}$ stands for the Bessel function of the first kind of order
zero. All these properties follow from standard Sturm-Liouville theory and
from related properties of Bessel functions (it is worth recalling here that
(\ref{neumann}) is Neumann's boundary condition at $r=1$ for the problem
under consideration since $J_{0}^{\prime }=-J_{1}$, see, e.g., Section 40 in
Chapter VII of \cite{weinberger}, and that the factor two in (\ref{neumann})
and (\ref{eigenfunctions}) is due to the factor one-half in (\ref%
{forwardproblembis})-(\ref{backwardproblembis})). For every $\mathsf{m}\in 
\mathbb{N}$ let us now choose the initial-final data as%
\begin{equation}
\varphi _{0,\mathsf{m}}(r)=\left\{ 
\begin{array}{c}
\frac{1}{\pi }\text{ \ for }\mathsf{m=0,} \\ 
\\ 
\frac{1}{\pi }\left( 1+J_{0}\left( \sqrt{2\lambda _{\mathsf{m}}}r\right)
\right) \text{ \ for }\mathsf{m}\in \mathbb{N}^{+}%
\end{array}%
\right.   \label{initialdata}
\end{equation}%
and%
\begin{equation}
\psi _{T,\mathsf{m}}(r)=1,  \label{finaldata}
\end{equation}%
respectively. It follows from (\ref{initialdata}) and (\ref{finaldata}) that
Hypothesis (H$_{2}$) holds, a consequence of elementary properties of $J_{0}$
including its uniform boundedness in case of (\ref{initialdata}). The
corresponding solutions to (\ref{forwardproblembis}) and (\ref%
{backwardproblembis}) then read

\bigskip

\begin{equation}
u_{\mathsf{m}}(r,t)=\left\{ 
\begin{array}{c}
\frac{1}{\pi }\text{ \ for }\mathsf{m=0,} \\ 
\\ 
\frac{1}{\pi }\left( 1+\exp \left[ -t\lambda _{\mathsf{m}}\right]
J_{0}\left( \sqrt{2\lambda _{\mathsf{m}}}r\right) \right) \text{ \ for }%
\mathsf{m}\in \mathbb{N}^{+}%
\end{array}%
\right.  \label{forwardsolutionbis}
\end{equation}%
and 
\begin{equation}
v_{\mathsf{m}}(r,t)=1  \label{backwardsolutionbis}
\end{equation}%
for each $r\in \left[ 0,1\right] $ and every $t\in \left[ 0,T\right] $,
respectively. Moreover, as a consequence of (\ref{backwardsolution}), (\ref%
{initialdata}) and (\ref{backwardsolutionbis}) we also have%
\begin{eqnarray*}
&&\int_{\left\vert \mathsf{x}\right\vert <1}\mathsf{dx}\varphi _{\mathsf{m,0}%
}(\mathsf{x)}\dint\limits_{\left\vert \mathsf{x}\right\vert <1}\mathsf{dy}g(%
\mathsf{x},T,\mathsf{y})\psi _{\mathsf{m,}T}(\mathsf{y}) \\
&=&\int_{\left\vert \mathsf{x}\right\vert <1}\mathsf{dx}\varphi _{\mathsf{m,0%
}}(\mathsf{x)=2\pi }\int_{0}^{1}drr\varphi _{\mathsf{m,0}}(r\mathsf{)=1}
\end{eqnarray*}%
so that (\ref{normalization}) is verified. We may therefore apply all the
results of the preceding section to the present situation, some of which we
state in the following proposition where%
\begin{equation*}
\mathbb{D}:=\left\{ \mathsf{x}\in \mathbb{R}^{2}:\left\vert \mathsf{x}%
\right\vert <1\right\} .
\end{equation*}

\bigskip

\textbf{Proposition 2.} \textit{For each }$\mathsf{m}\in \mathbb{N}$\textit{%
\ let }$\mu _{\mathsf{m}}$\textit{\ be the measure of the form (\ref%
{probabilitymeasures}) with the initial-final data given by (\ref%
{initialdata}) and (\ref{finaldata}), respectively. Then, there exists a }$%
\overline{\mathbb{D}}$\textit{-valued, non-stationary Markovian Bernstein
process} $Z_{\tau \in \left[ 0,T\right] }^{\mathsf{m}}$ \textit{such that
the following properties hold:}

\textit{(a) For each Borel subset }$F\subseteq \overline{\mathbb{D}}$ 
\textit{of Lebesgue measure }$\left\vert F\right\vert $\textit{\ and for
every }$t\in \left[ 0,T\right] $\textit{\ we have}%
\begin{equation*}
\mathbb{P}_{\mu _{\mathsf{m}}}\left( Z_{t}^{\mathsf{m}}\in F\right) =\left\{ 
\begin{array}{c}
\frac{\left\vert F\right\vert }{\pi }\text{ \ for }\mathsf{m=0,} \\ 
\\ 
\frac{1}{\pi }\left( \left\vert F\right\vert +\exp \left[ -t\lambda _{%
\mathsf{m}}\right] \int_{F}\mathsf{dx}J_{0}\left( \sqrt{2\lambda _{\mathsf{m}%
}}\left\vert \mathsf{x}\right\vert \right) \right) \text{ \ for }\mathsf{m}%
\in \mathbb{N}^{+}.%
\end{array}%
\right.
\end{equation*}%
\textit{Thus, the function} $t\rightarrow \mathbb{P}_{\mu _{\mathsf{m}%
}}\left( Z_{t}^{\mathsf{m}}\in F\right) $ \textit{is non-increasing on} $%
\left[ 0,T\right] $.

\textit{(b)} \textit{For each bounded Borel measurable function }$b:%
\overline{\mathbb{D}}$ $\mathbb{\mapsto C}$ \textit{and every }$t\in \left[
0,T\right] $ \textit{we have}%
\begin{equation*}
\mathbb{E}_{\mu _{\mathsf{m}}}\left( b(Z_{t}^{\mathsf{m}})\right) =\left\{ 
\begin{array}{c}
\frac{1}{\pi }\int_{\left\vert \mathsf{x}\right\vert <1}\mathsf{dx}b(\mathsf{%
x})\text{ \ for }\mathsf{m=0,} \\ 
\\ 
\frac{1}{\pi }\int_{\left\vert \mathsf{x}\right\vert <1}\mathsf{dx}b(\mathsf{%
x})\left( 1+\exp \left[ -t\lambda _{\mathsf{m}}\right] J_{0}\left( \sqrt{%
2\lambda _{\mathsf{m}}}\left\vert \mathsf{x}\right\vert \right) \right) 
\text{ \ for }\mathsf{m}\in \mathbb{N}^{+}.%
\end{array}%
\right.
\end{equation*}

\textit{(c) Let }$\lambda \in \mathbb{R}^{+}$ \textit{be given. Then, the
process }$\bar{Z}_{\tau \in \left[ 0,T\right] }$ \textit{of maximal entropy
within }$\overline{\mathbb{D}}$ \textit{in the sense of Theorem 4} \textit{is%
} \textit{obtained by averaging the }$Z_{\tau \in \left[ 0,T\right] }^{%
\mathsf{m}}$\textit{\ with probabilities of the form}%
\begin{equation*}
p_{\mathsf{m}}=\mathcal{Z}^{-1}\mathcal{(}\beta \left( \lambda \right) )\exp %
\left[ -\beta \left( \lambda \right) \lambda _{\mathsf{m}}\right]
\end{equation*}%
\textit{where }$\mathcal{Z(}\beta \left( \lambda \right) )$\textit{\ is
given by (\ref{partitionfunctionbis}). Moreover, }$\bar{Z}_{\tau \in \left[
0,T\right] }$\textit{\ is Markovian and its entropy may be evaluated from (%
\ref{maxentropy}).}

\bigskip

The proof is a direct application of the corresponding formulae in Section 2
combined with those of this section. We note that we must have $\lambda >0$
for the preassigned value in Statement (c) since $\lambda _{\mathsf{m}}\geq
0 $ for every $\mathsf{m}\in \mathbb{N}$ according to (\ref{ordering}). We
also have $\sum_{\mathsf{m}\in \mathbb{N}^{d}}\exp \left[ -\beta \lambda _{%
\mathsf{m}}\right] \lambda _{\mathsf{m}}<+\infty $ for every $\beta \in
\left( 0,+\infty \right) $ as a consequence of (\ref{growth}), so that
expression (\ref{maxentropy}) may indeed be applied in this case.

\bigskip

\textsc{Remarks.} (1) Although the mixed processes obtained by the averaging
method described in Section 2 are not Markovian in general (see the remark
following the proof of Theorem 2), the reason why $\bar{Z}_{\tau \in \left[
0,T\right] }$ possesses the Markov property is due to our choice of the
final condition (\ref{finaldata}), which implies that the averaged joint
density (\ref{jointdensity}) becomes%
\begin{equation*}
\overline{\mu }\left( \mathsf{x,y}\right) =g(\mathsf{x},T,\mathsf{y})\sum_{%
\mathsf{m\in }\mathbb{N}}p_{\mathsf{m}}\varphi _{\mathsf{m,0}}\left( \mathsf{%
x}\right)
\end{equation*}%
where $\varphi _{\mathsf{m,0}}$ is given by (\ref{initialdata}). Indeed, the
preceding relation is then of the form (\ref{markovmeasure}) with an obvious
choice for $\nu _{0}$ and $\nu _{1}$. For more details and examples
regarding the time evolution of Bernstein processes that possess the Markov
property we refer the reader to \cite{vuillermot}.

(2) We can obtain similar results for radially symmetric forward-backward
problems of the form (\ref{forwardproblem})-(\ref{backwardproblem}) with $%
V=0 $ defined in the open ball \textsf{B}$^{d}\subset \mathbb{R}^{d}$ of
radius one centered at the origin where $d\geq 3$. The eigenfunctions of the
radial part of the Laplacian then involve the Bessel function $J_{\frac{d}{2}%
-1}$ rather than $J_{0}$, while Neumann's boundary condition is expressed in
terms of $J_{\frac{d}{2}}$ instead of $J_{1}$ (the case $d=1$ can be dealt
with directly in terms of trigonometric functions). However, the
corresponding formulae for the probabilities and the expectation values of
the underlying processes become much more involved.

\bigskip

\textbf{Acknowledgements.} The author would like to thank the Funda\c{c}\~{a}%
o para a Ci\^{e}ncia e Tecnologia (FCT) of the Portuguese Government for its
financial support under Grant PTDC/MAT-STA/0975/2014. Parts of this work
were presented as an invited contributed talk at the 2nd Conference of the
European Physical Society Statistical and Nonlinear Physics Division, which
took place at Nordita, Stockholm, in May of 2019. The author would like to
thank the organizers for their very kind invitation.

\bigskip

\textbf{Appendix A. A class of time-dependent statistical operators}

In this appendix we define and investigate statistical operators which are
more general than that defined by (\ref{statisticaloperator}), in view of
getting expressions such as (\ref{stataverageter}) for the averaged
processes of Theorem 2 which are as a rule non-stationary. This, however,
requires some additional structure. Let us write%
\begin{equation}
\exp \left[ -t\mathcal{H}\right] f\left( .\right) :=\left\{ 
\begin{array}{c}
f(.)\text{ \ \ if }t=0, \\ 
\\ 
\int_{D}\mathsf{dy}g(\mathsf{.},t,\mathsf{y})f(\mathsf{y)}\text{ \ \ if }%
t\in \left( 0,T\right]%
\end{array}%
\right.  \label{semigroup}
\end{equation}%
for the positivity preserving, symmetric semigroup generated by (\ref%
{hamiltonian}) on $L^{2}(D)$, where $g$ is given by (\ref{heatkernel}). In
addition to Hypothesis (H$_{2}$) regarding the initial-final data we shall
now impose the following requirement:

\bigskip

(H$_{3}$) The sequences $\left( \varphi _{\mathsf{m,0}}\right) _{\mathsf{%
m\in }\mathbb{N}^{d}}$ and $\left( \exp \left[ -T\mathcal{H}\right] \psi _{%
\mathsf{m,}T}\right) _{\mathsf{m\in }\mathbb{N}^{d}}$ form a biorthonormal
system in $L^{2}(D)$, that is,%
\begin{equation*}
\left( \varphi _{\mathsf{m,0}},\exp \left[ -T\mathcal{H}\right] \psi _{%
\mathsf{n,}T}\right) _{2}=\delta _{\mathsf{m},\mathsf{n}}
\end{equation*}%
for all\textsf{\ }$\mathsf{m,n}\in \mathbb{N}^{d}$.

\bigskip First we have:

\textbf{Lemma 2.}\textit{\ Let us assume that Hypotheses (}H$_{2}$\textit{)
and (}H$_{3}$\textit{) hold. Then the sequences of solutions }$\left( u_{%
\mathsf{m}}(\mathsf{.},t)\right) _{\mathsf{m\in }\mathbb{N}^{d}}$,\textit{\ }%
$\left( v_{\mathsf{m}}(\mathsf{.},t)\right) _{\mathsf{m\in }\mathbb{N}^{d}}$%
\textit{\ given by (\ref{forwardsolution}) and (\ref{backwardsolution}),
respectively, form a biorthonormal system in }$L^{2}(D)$\textit{\ for every }%
$t\in \left[ 0,T\right] $\textit{.}

\bigskip

\textbf{Proof.} From (\ref{forwardsolution}), (\ref{backwardsolution}) and (%
\ref{semigroup}) we have%
\begin{eqnarray}
&&\left( u_{\mathsf{m}}(\mathsf{.},t),v_{\mathsf{n}}(\mathsf{.},t)\right)
_{2}  \notag \\
&=&\left( \exp \left[ -t\mathcal{H}\right] \varphi _{\mathsf{m,0}},\exp %
\left[ -(T-t)\mathcal{H}\right] \psi _{\mathsf{n,}T}\right) _{2}
\label{biorthogonality} \\
&=&\left( \varphi _{\mathsf{m,0}},\exp \left[ -T\mathcal{H}\right] \psi _{%
\mathsf{n,}T}\right) _{2}=\delta _{\mathsf{m},\mathsf{n}}  \notag
\end{eqnarray}%
from the symmetry of the semigroup and (H$_{3}$), for all $\mathsf{m},%
\mathsf{n}\in \mathbb{N}^{d}$. \ \ $\blacksquare $\ 

\bigskip

We then define%
\begin{equation}
\mathcal{R}\left( t\right) f:=\sum_{\mathsf{m\in }\mathbb{N}^{d}}p_{\mathsf{m%
}}\left( f,u_{\mathsf{m}}(\mathsf{.},t)\right) _{2}v_{\mathsf{m}}(\mathsf{.}%
,t)  \label{statisticaloperatorbis}
\end{equation}%
for each $f\in L^{2}(D)$ and every $t\in \left[ 0,T\right] $. The following
result shows that (\ref{statisticaloperatorbis}) possesses several
properties similar to those stated in Proposition 1:

\bigskip

\textbf{Proposition A.1.} \textit{Let us assume that Hypotheses (}H$_{1}$%
\textit{), (}H$_{2}$\textit{) and (}H$_{3}$\textit{) hold. Then the
following statements are valid:}

\textit{(a) Expression (\ref{statisticaloperatorbis}) defines a linear
trace-class operator in }$L^{2}\left( D\right) $ \textit{such that}%
\begin{equation}
\func{Tr}\mathcal{R}=\sum_{\mathsf{m\in }\mathbb{N}^{d}}p_{\mathsf{m}}=1
\label{traceter}
\end{equation}%
\textit{and}%
\begin{equation}
\func{Tr}\mathcal{R}^{2}=\sum_{\mathsf{m\in }\mathbb{N}^{d}}p_{\mathsf{m}%
}^{2}\leq 1  \label{tracequarto}
\end{equation}%
\textit{for every} $t\in \left[ 0,T\right] $.

\bigskip \textit{(b) The eigenvalue equation}%
\begin{equation}
\mathcal{R}\left( t\right) v_{\mathsf{m}}(\mathsf{.},t)=p_{\mathsf{m}}v_{%
\mathsf{m}}(\mathsf{.},t)  \label{eigenvalueequationbis}
\end{equation}%
\textit{holds for each} $\mathsf{m}\in \mathbb{N}^{d}$ \textit{and every} $%
t\in \left[ 0,T\right] $.

\textit{(c) If }$B$\textit{\ is a linear bounded self-adjoint operator on }$%
L^{2}(D)$ \textit{we have}%
\begin{equation}
\func{Tr}\left( \mathcal{R}(t)B\right) =\sum_{\mathsf{m\in }\mathbb{N}%
^{d}}p_{\mathsf{m}}\left( Bu_{\mathsf{m}}(\mathsf{.},t),v_{\mathsf{m}}(%
\mathsf{.},t)\right) _{2}  \label{staaveragequinto}
\end{equation}%
\textit{for every} $t\in \left[ 0,T\right] $. \textit{In particular, if }$B$%
\textit{\ is the same multiplication operator as in Part (c) of Proposition
1 we have}%
\begin{equation}
\func{Tr}\left( \mathcal{R}(t)B\right) \mathcal{=}\mathbb{E}_{_{\bar{\mu}%
}}\left( b(\bar{Z}_{t})\right)  \label{stataveragequarto}
\end{equation}%
\textit{where the right-hand side of (\ref{stataveragequarto}) is given by (%
\ref{expectationbis}).}

\bigskip

\textbf{Proof.} The proof of the trace-class property is similar to that
given in Proposition 1. Thus, we first remark that there exists a finite
constant $c>0$ such that the estimates%
\begin{equation*}
\left\Vert u_{\mathsf{m}}(\mathsf{.},t)\right\Vert _{2}\leq c
\end{equation*}%
and%
\begin{equation*}
\left\Vert v_{\mathsf{m}}(\mathsf{.},t)\right\Vert _{2}\leq c
\end{equation*}%
hold uniformly in $\mathsf{m}$ and $t$ as a consequence of (\ref{heatkernel}%
), (\ref{forwardsolution}), (\ref{backwardsolution}) and Hypothesis (H$_{2}$%
), which makes a linear bounded operator out of (\ref{statisticaloperatorbis}%
) on $L^{2}(D)$. The relevant auxiliary function for the remaining part of
the argument is then given by%
\begin{equation}
a\left( \mathsf{m,n},t\right) :=p_{\mathsf{m}}\left( \mathsf{h}_{\mathsf{n}%
},u_{\mathsf{m}}(\mathsf{.},t)\right) _{2}\left( v_{\mathsf{m}}(\mathsf{.}%
,t),\mathsf{h}_{\mathsf{n}}\right) _{2}  \label{auxiliaryfunctionbis}
\end{equation}%
for every $t\in \left[ 0,T\right] $, with $\left( \mathsf{h}_{\mathsf{n}%
}\right) _{\mathsf{n\in }\mathbb{N}^{d}}$ an arbitrary orthogonal basis as
before. It is indeed easily seen that the properties of (\ref%
{auxiliaryfunctionbis}) are similar to those of (\ref{auxiliaryfunction}),
the key point in getting (\ref{traceter}) and (\ref{tracequarto}) being the
biorthogonality relation (\ref{biorthogonality}) which replaces the
orthogonality properties of $\left( \mathsf{f}_{\mathsf{m}}\right) _{\mathsf{%
m}\in \mathbb{N}^{d}}$. The same observation applies for the proof of (\ref%
{eigenvalueequationbis}), while (\ref{staaveragequinto}) follows from the
relation%
\begin{equation*}
\left( \mathcal{R}(t)B\mathsf{h}_{\mathsf{n}},\mathsf{h}_{\mathsf{n}}\right)
_{2}=\sum_{\mathsf{m\in }\mathbb{N}^{d}}p_{\mathsf{m}}\left( \mathsf{h}_{%
\mathsf{n}},Bu_{\mathsf{m}}(\mathsf{.},t)\right) _{2}(v_{\mathsf{m}}(\mathsf{%
.},t),\mathsf{h}_{\mathsf{n}})_{2}
\end{equation*}%
valid for every $\mathsf{n}\in \mathbb{N}^{d}$, Relation (\ref%
{staaveragequinto}) then implying (\ref{stataveragequarto}). \ \ $%
\blacksquare $

\bigskip

\textsc{Remarks.} (1) As a consequence of (\ref{semigroup}) and elementary
spectral theory, it is plain that the forward and backward solutions (\ref%
{forwardsolution}) and (\ref{backwardsolution}) become%
\begin{equation*}
u_{\mathsf{m}}(\mathsf{.},t)=\exp \left[ -t\lambda _{\mathsf{m}}\right] 
\mathsf{f}_{\mathsf{m}}
\end{equation*}%
and 
\begin{equation*}
v_{\mathsf{m}}(\mathsf{.},t)=\exp \left[ t\lambda _{\mathsf{m}}\right] 
\mathsf{f}_{\mathsf{m}},
\end{equation*}%
respectively, for initial-final data given by (\ref{initialfinaldata}). The
substitution of these expressions into (\ref{statisticaloperatorbis}) then
gives (\ref{statisticaloperator}), so that the former operators are indeed
time-dependent generalizations of the latter. Moreover, thanks to Statement
(a) of Proposition A.1 the classification of processes as pure or mixed
according to Definition 2 still holds for the more general form (\ref%
{statisticaloperatorbis}).

(2) We may write (\ref{statisticaloperatorbis}) as%
\begin{equation*}
\mathcal{R}\left( t\right) f=\sum_{\mathsf{m\in }\mathbb{N}^{d}}p_{\mathsf{m}%
}\mathcal{P}_{\mathsf{m}}\left( t\right) f
\end{equation*}%
where the operators given by%
\begin{equation*}
\mathcal{P}_{\mathsf{m}}\left( t\right) f:=\left( f,u_{\mathsf{m}}(\mathsf{.}%
,t)\right) _{2}v_{\mathsf{m}}(\mathsf{.},t)
\end{equation*}%
satisfy%
\begin{equation*}
\mathcal{P}_{\mathsf{m}}^{2}\left( t\right) =\mathcal{P}_{\mathsf{m}}\left(
t\right)
\end{equation*}%
for each $\mathsf{m}\in \mathbb{N}^{d}$ and every $t\in \left[ 0,T\right] $
as a consequence of the biorthogonality property. Therefore, we may view $%
\mathcal{R}\left( t\right) $ as a statistical mixture of oblique projections
as the $\mathcal{P}_{\mathsf{m}}\left( t\right) $ are not self-adjoint in
the general case. In fact, the adjoint of (\ref{statisticaloperatorbis}) in $%
L^{2}(D)$ is obtained by swapping the r\^{o}les of $u_{\mathsf{m}}(\mathsf{.}%
,t)$ and $v_{\mathsf{m}}(\mathsf{.},t)$, that is, 
\begin{equation}
\mathcal{R}^{\ast }\left( t\right) f:=\sum_{\mathsf{m\in }\mathbb{N}^{d}}p_{%
\mathsf{m}}\left( f,v_{\mathsf{m}}(\mathsf{.},t)\right) _{2}u_{\mathsf{m}}(%
\mathsf{.},t).  \label{adjoint}
\end{equation}%
We note that (\ref{adjoint}) enjoys the very same properties as (\ref%
{statisticaloperatorbis}), with the exception of (\ref{eigenvalueequationbis}%
) which has to be replaced by%
\begin{equation*}
\mathcal{R}^{\ast }\left( t\right) u_{\mathsf{m}}(\mathsf{.},t)=p_{\mathsf{m}%
}u_{\mathsf{m}}(\mathsf{.},t).
\end{equation*}%
We could therefore have chosen (\ref{adjoint}) as statistical operators
instead of (\ref{statisticaloperatorbis}). Finally, we remark that $\mathcal{%
R}\left( t\right) $ and $\mathcal{R}^{\ast }\left( t\right) $ both involve (%
\ref{forwardsolution}) and (\ref{backwardsolution}), in agreement with the
fact that there are two time directions in the theory from the outset.

(3) It is reasonable to ask whether it is always possible to choose
initial-final data so that Hypothesis (H$_{3}$) holds. The answer is
affirmative provided the functions $\exp \left[ -T\mathcal{H}\right] \psi _{%
\mathsf{m,}T}$ remain close to the orthonormal basis $\left( \mathsf{f}_{%
\mathsf{m}}\right) _{\mathsf{m\in }\mathbb{N}^{d}}$ in some very specific $%
L^{2}(D)$-sense. This follows from a direct application of the
generalization of a theorem by Paley and Wiener as stated in Section 86 of
Chapter V in \cite{riesznagy}. In that case the sequences $\left( \varphi _{%
\mathsf{m,0}}\right) _{\mathsf{m\in }\mathbb{N}^{d}}$ and $\left( \exp \left[
-T\mathcal{H}\right] \psi _{\mathsf{m,}T}\right) _{\mathsf{m\in }\mathbb{N}%
^{d}}$ form a \textit{complete} biorthogonal system.

\bigskip

\textbf{Appendix B. On the existence of Bernstein processes and their
relation with Schr\"{o}dinger's problem and Optimal Transport Theory}

The typical construction of a Bernstein process with state space $\overline{D%
}$ requires a probability measure $\mu $ on $\mathcal{B}(\overline{D})\times 
\mathcal{B}(\overline{D})$ and a transition function $Q$, as is the case for
Markov processes. We provide below a general theorem which is a direct
consequence of a more abstract construction carried out in \cite{jamison},
or with a more analytical flavor in \cite{vuillerzambis}, to which we refer
the reader for details. The theorem shows that all the basic quantities that
characterize a Bernstein process can be expressed exclusively in terms of $%
\mu $ and the heat kernel $g$, which is all we needed in the preceding
sections.

Since there are two time directions provided by (\ref{forwardproblem})-(\ref%
{backwardproblem}), the natural choice for the transition function we
alluded to is%
\begin{equation}
Q\left( \mathsf{x},t;F,r;\mathsf{y},s\right) :=\dint\limits_{F}\mathsf{dz}%
q\left( \mathsf{x},t;\mathsf{z},r;\mathsf{y},s\right)
\label{transitionfunction}
\end{equation}%
for every $F\in \mathcal{B}(\overline{D})$, where%
\begin{equation}
q\left( \mathsf{x},t;\mathsf{z},r;\mathsf{y},s\right) :=\frac{g(\mathsf{x}%
,t-r,\mathsf{z})g(\mathsf{z},r-s,\mathsf{y})}{g(\mathsf{x},t-s,\mathsf{y})}.
\end{equation}%
Both functions are well defined and positive for all $\mathsf{x},\mathsf{y},%
\mathsf{z}\in \mathbb{R}^{d}$ and all $r,s,t$ satisfying $r\in \left(
s,t\right) \subset \left[ 0,T\right] $ by virtue of (\ref{heatkernel}), and
moreover the normalization condition%
\begin{equation*}
Q\left( \mathsf{x},t;D,r;\mathsf{y},s\right) =1
\end{equation*}%
holds as a consequence of the semigroup composition law for $g$. The result
we have in mind is the following:

\bigskip

\textbf{Theorem B.1.} \textit{Let }$\mu $\textit{\ be a probability measure
on }$\mathcal{B}(\overline{D})\times \mathcal{B}(\overline{D})$, \textit{and
let }$Q$ \textit{be given by (\ref{transitionfunction}). Then there exists a
probability space }$\left( \Omega ,\mathcal{F},\mathbb{P}_{\mu }\right) $ 
\textit{supporting an }$\overline{D}$\textit{-valued Bernstein process }$%
Z_{\tau \in \left[ 0,T\right] }$\textit{\ such that the following properties
are valid:}

\textit{(a) The function }$Q$\textit{\ is the two-sided transition function
of }$Z_{\tau \in \left[ 0,T\right] }$\textit{\ in the sense that} 
\begin{equation}
\mathbb{P}_{\mu }\left( Z_{r}\in F\left\vert Z_{s},Z_{t}\right. \right)
=Q\left( Z_{t},t;F,r;Z_{s},s\right)
\end{equation}%
\textit{for each }$F\in \mathcal{B}(\overline{D})$ \textit{and all }$r,s,t$%
\textit{\ satisfying }$r\in \left( s,t\right) \subset \left[ 0,T\right] $. 
\textit{Moreover,}%
\begin{equation}
\mathbb{P}_{\mu }\left( Z_{0}\in F_{0},Z_{T}\in F_{T}\right) =\mu \left(
F_{0}\times F_{T}\right)
\end{equation}%
\textit{for all }$F_{0},F_{T}\in \mathcal{B}(\overline{D})$\textit{, that
is, }$\mu $\textit{\ is the joint probability distribution of }$Z_{0}$%
\textit{\ and }$Z_{T}$\textit{. In particular, the marginal distributions
are given by}%
\begin{equation}
\mathbb{P}_{\mu }\left( Z_{0}\in F\right) =\mu \left( F\times \overline{D}%
\right)
\end{equation}%
\textit{and}%
\begin{equation}
\mathbb{P}_{\mu }\left( Z_{T}\in F\right) =\mu \left( \overline{D}\times
F\right)
\end{equation}%
\textit{for each} $F\in \mathcal{B}(\overline{D})$, \textit{respectively.}

\textit{(b) For every }$n\in \mathbb{N}^{+}$ \textit{the finite-dimensional
distributions of the process are given by}%
\begin{eqnarray}
&&\mathbb{P}_{\mu }\left( Z_{t_{1}}\in F_{1},...,Z_{t_{n}}\in F_{n}\right) 
\notag \\
&=&\int_{D\times D}\frac{\mathsf{d}\mu \mathsf{\left( \mathsf{x,y}\right) }}{%
g(\mathsf{x},T,\mathsf{y})}\int_{F_{1}}\mathsf{dx}_{1}...\int_{F_{n}}\mathsf{%
dx}_{n}  \notag \\
&&\times \dprod\limits_{k=1}^{n}g\left( \mathsf{x}_{k},t_{k}-t_{k-1},\mathsf{%
x}_{k-1}\right) \times g\left( \mathsf{y},T-t_{n},\mathsf{x}_{n}\right)
\end{eqnarray}%
\textit{for all }$F_{1},...,F_{n}\in $\textit{\ }$\mathcal{B}(\overline{D})$ 
\textit{and all }$t_{0}=0<t_{1}<...<t_{n}<T$\textit{, where }$\mathsf{x}_{0}=%
\mathsf{x}$\textit{. In particular we have}%
\begin{eqnarray}
&&\mathbb{P}_{\mu }\left( Z_{t}\in F\right)  \notag \\
&=&\int_{D\times D}\frac{\mathsf{d}\mu \mathsf{\left( \mathsf{x,y}\right) }}{%
g(\mathsf{x},T,\mathsf{y})}\int_{F}\mathsf{dz}g\left( \mathsf{x},t,\mathsf{z}%
\right) g\left( \mathsf{z},T-t,\mathsf{y}\right)
\end{eqnarray}%
\textit{for each }$F\in \mathcal{B}(\overline{D})$\textit{\ and every} $t\in
\left( 0,T\right) $.

\textit{(c) }$\mathbb{P}_{\mu }$\textit{\ is the only probability measure
leading to the above properties.}

\bigskip

As we saw in Section 2, Theorems 1, 2 and 3 were obtained by substituting
the respective measures in the formulae of Theorem B.1.

\bigskip

We conclude this appendix with a very brief remark which establishes the
connection between Markovian Bernstein processes on the one hand, Schr\"{o}%
dinger's problem and Optimal Transport Theory on the other hand. From
Section 2 we know that in the Markovian case the relevant probability
measures to be substituted into the formulae of Theorem B.1 are necessarily
of the form%
\begin{equation*}
\mu (G)=\int_{G}\mathsf{dxdy}\varphi _{\mathsf{0}}(\mathsf{x)}g(\mathsf{x},T,%
\mathsf{y})\psi _{T}(\mathsf{y})
\end{equation*}%
where $G\in \mathcal{B}(\overline{D})\times \mathcal{B}(\overline{D})$, and
where $\varphi _{\mathsf{0}}>0$ and $\psi _{T}>0$. In particular, the
marginal distributions are%
\begin{equation*}
\mu (F\times \overline{D})=\int_{F}\mathsf{dx}\varphi _{\mathsf{0}}(\mathsf{%
x)}\int_{D}\mathsf{dy}g(\mathsf{x},T,\mathsf{y})\psi _{T}(\mathsf{y})
\end{equation*}%
and 
\begin{equation*}
\mu (\overline{D}\times F)=\int_{F}\mathsf{dy}\psi _{T}(\mathsf{y})\int_{D}%
\mathsf{dx}\varphi _{\mathsf{0}}(\mathsf{x)}g(\mathsf{x},T,\mathsf{y})
\end{equation*}%
where $F\in \mathcal{B}(\overline{D})$, which gives rise to the respective
densities%
\begin{equation}
\mu _{0}(\mathsf{x}):=\varphi _{\mathsf{0}}(\mathsf{x)}\int_{D}\mathsf{dy}g(%
\mathsf{x},T,\mathsf{y})\psi _{T}(\mathsf{y})  \label{density1}
\end{equation}%
and 
\begin{equation}
\mu _{T}(\mathsf{y}):=\psi _{T}(\mathsf{y})\int_{D}\mathsf{dx}\varphi _{%
\mathsf{0}}(\mathsf{x)}g(\mathsf{x},T,\mathsf{y}).  \label{density2}
\end{equation}%
Thus, the considerations of this article show that these marginal densities
are entirely determined by the initial-final data once the heat kernel is
known. It is, however, the inverse point of view that prevails in Schr\"{o}%
dinger's problem and in the related Optimal Transport Theory, which \textit{%
first }amounts to prescribing $\mu _{0}$ and $\mu _{T}$ and \textit{then}
consider (\ref{density1}) and (\ref{density2}) as a nonlinear inhomogeneous
system of integral equations in the two unknowns $\varphi _{\mathsf{0}}$ and 
$\psi _{T}$. That is actually what was developed by Schr\"{o}dinger in the
last part of \cite{schroedinger} by using entropy arguments. Moreover, given 
$\mu _{0}$ and $\mu _{T}$ continuous, a very general existence and
uniqueness result for the pair $\left( \varphi _{\mathsf{0}},\psi
_{T}\right) $ satisfying (\ref{density1}) and (\ref{density2}) was proved in 
\cite{beurling}. In the context of Optimal Transport Theory, arguments that
allow the minimization of the so-called cost functions using entropy related
methods are often used. We refer the reader fro instance to \cite{leonard}
and to the references therein for details.

\end{document}